
\def\LaTeX{%
  \let\Begin\begin
  \let\End\end
  \let\salta\relax
  \let\finqui\relax
  \let\futuro\relax}

\def\UK{\def\our{our}\let\sz s}
\def\USA{\def\our{or}\let\sz z}

\LaTeX

\USA


\salta

\documentclass[twoside,12pt]{article}
\setlength{\textheight}{24cm}
\setlength{\textwidth}{16cm}
\setlength{\oddsidemargin}{2mm}
\setlength{\evensidemargin}{2mm}
\setlength{\topmargin}{-15mm}
\parskip2mm


\usepackage{amsmath}
\usepackage{amsthm}
\usepackage{amssymb}
\usepackage[mathcal]{euscript}
\usepackage{cite}

%
%
%
\def\gianni #1{{\color{blue}#1}}
\def\juerg #1{{\color{blue}#1}}          
\def\pier #1{{\color{red}#1}}     
%
%
\let\gianni\relax
\let\juerg\relax
\let\pier\relax




\bibliographystyle{plain}


%

\finqui

\newcommand{\beq}{\begin{equation}}
\newcommand{\eeq}{\end{equation}}
\newcommand{\beqa}{\begin{eqnarray}}
\newcommand{\eeqa}{\end{eqnarray}}
\let\non\nonumber


\def\step #1 \par{\medskip\noindent{\bf #1.}\quad}


\def\aand{\quad\hbox{and}\quad}

\def\lhs{left-hand side}
\def\rhs{right-hand side}


\def\multibold #1{\def\arg{#1}%
  \ifx\arg\pto \let\next\relax
  \else
  \def\next{\expandafter
    \def\csname #1#1#1\endcsname{{\bf #1}}%
    \multibold}%
  \fi \next}

\def\pto{.}

\def\multical #1{\def\arg{#1}%
  \ifx\arg\pto \let\next\relax
  \else
  \def\next{\expandafter
    \def\csname cal#1\endcsname{{\cal #1}}%
    \multical}%
  \fi \next}


\def\multimathop #1 {\def\arg{#1}%
  \ifx\arg\pto \let\next\relax
  \else
  \def\next{\expandafter
    \def\csname #1\endcsname{\mathop{\rm #1}\nolimits}%
    \multimathop}%
  \fi \next}

\multibold
qwertyuiopasdfghjklzxcvbnmQWERTYUIOPASDFGHJKLZXCVBNM.

\multical
QWERTYUIOPASDFGHJKLZXCVBNM.

\multimathop
dist div dom meas sign supp .


\def\Accorpa #1#2 #3 {\gdef #1{\eqref{#2}--\eqref{#3}}%
  \wlog{}\wlog{\string #1 -> #2 - #3}\wlog{}}


\def\qed{\hfill $\square$}
\def\tint {\int_0^t}

\def\xinto{\int_\Omega}
\def\txinto{\int_0^t\!\!\!\int_\Omega}

\def\oma{\Omega}
\def\ginto{\int_\Gamma}
\def\tginto{\int_0^t\!\!\!\int_\Gamma}

\def\pt{\partial_t}
\def\pn{\partial_{\bf n}}
\def\ds{\,{\rm d}s}
\def\dt{\,{\rm d}t}
\def\dx{\,{\rm d}x}

\def\ug{\/ u_\Gamma}

\def\mun{\mu_n}
\def\rhon{\rho_n}
\def\rhong{\rho_{n_\Gamma}}

\def\munmu{\mu_{n-1}}

\def\rhone{\rho_n^\eps}

\def\checkmmode #1{\relax\ifmmode\hbox{#1}\else{#1}\fi}

\def\aeQ{\checkmmode{a.\,e.\ in~$Q$}}

\def\limrho0{\lim_{\rho\searrow 0}}
\def\limrho1{\lim_{\rho\nearrow 1}}

\def\rz{{\mathbb{R}}}
\def\nz{{\mathbb{N}}}


\def\Lx #1{L^{#1}(\Omega)}
\def\Hx #1{H^{#1}(\Omega)}

\def\Ldue{\Lx 2}

\def\Huno{\Hx 1}
\def\Hdue{\Hx 2}

\def\Hg{H_\Gamma}
\def\Vg{V_\Gamma}



\let\eps\varepsilon

\let\TeXchi\chi                         
\newbox\chibox
\setbox0 \hbox{\mathsurround0pt $\TeXchi$}
\setbox\chibox \hbox{\raise\dp0 \box 0 }
\def\chi{\copy\chibox}


\def\r0g{\rho_{0_{|\Gamma}}}

\def\rg{\rho_\Gamma}

\def\pig{\pi_\Gamma}
\def\xig{\xi_\Gamma}
\def\beg{\beta_\Gamma}
\def\ng{\nabla_\Gamma}
\def\Dg{\Delta_\Gamma}
\def\dg{\,{\rm d}\Gamma}

\def\coeff{1+2g(\rho)}
\def\coeffs{1+2g(\rho(s))}

\def\T{\calT_\tau}

\def\mute{\mu^\eps_\tau}
\def\rhote{\rho^\eps_\tau}
\def\rhotge{\rho^\eps_{\tau_\Gamma}}

\def\mun{\mu_n}
\def\rhon{\rho_n}

\def\munmu{\mu_{n-1}}

\Begin{document}


\title{\bf Global existence for a nonstandard viscous
Cahn--Hilliard system with dynamic boundary condition\footnote{This work received a partial support from the GNAMPA (Gruppo Nazionale per l'Analisi
Matematica, la Probabilit\`{a} e loro Applicazioni) of INDAM (Istituto Nazionale
di Alta Matematica) and the \pier{IMATI -- C.N.R. Pavia} for PC and GG.}}

\author{}
\date{}
\maketitle
\begin{center}
\vskip-2cm
{\large\bf Pierluigi Colli$^{(1)}$}\\
{\normalsize e-mail: {\tt pierluigi.colli@unipv.it}}\\[.25cm]
{\large\bf Gianni Gilardi$^{(1)}$}\\
{\normalsize e-mail: {\tt gianni.gilardi@unipv.it}}\\[.25cm]
{\large\bf J\"urgen Sprekels$^{(2)}$}\\
{\normalsize e-mail: {\tt sprekels@wias-berlin.de}}\\[.45cm]
$^{(1)}$
{\small Dipartimento di Matematica ``F. Casorati'', Universit\`a di Pavia}\\
{\small via Ferrata \gianni 5, 27100 Pavia, Italy}\\[.2cm]
$^{(2)}$
{\small Weierstrass Institute for Applied Analysis and Stochastics}\\
{\small Mohrenstra\ss e\ 39, 10117 Berlin, Germany, and}\\
{\small Department of Mathematics, Humboldt-Universit\"at zu Berlin}\\
{\small Unter den Linden 6, 10099 Berlin. Germany}\\[.8cm]
\end{center}


\Begin{abstract}\noindent
In this paper, we study a model for phase \pier{segregation} 
taking place in a spatial domain that was introduced 
\gianni{by Podio-Guidugli in Ric.\ Mat.\ {\bf 55} (2006), pp.~105--118}.
The model consists of a strongly  coupled system of nonlinear
parabolic differential equations, in  which products between  
the unknown functions and their time derivatives occur that are difficult to handle
analytically. In contrast to the existing literature about this PDE system, we 
consider here a dynamic boundary condition involving
the Laplace--Beltrami operator for the order parameter\pier{. This 
boundary condition} models an additional nonconserving phase
transition occurring on the surface of the domain. Different well-posedness
results are shown, depending on the smoothness properties of the involved
bulk and surface free energies. 
\\[7mm]
\End{abstract}


\salta

\pagestyle{myheadings}
\newcommand\testopari{\sc Colli \ --- \ Gilardi \ --- \ Sprekels}
\newcommand\testodispari{\sc\small Nonstandard
Cahn--Hilliard system with dynamic boundary condition}
\markboth{\testodispari}{\testopari}

\finqui


\section{Introduction}
\label{Intro}
\setcounter{equation}{0}
Let $\oma\subset\rz^3$ be a \gianni{bounded and connected open set} with a smooth boundary $\Gamma$, as well 
as $Q:=\oma\times (0,T)$ and $\Sigma:=\Gamma\times (0,T)$. We denote by $\pn$, $\ng$,
$\Delta_\Gamma$, the outward normal derivative, the tangential gradient, and the Laplace--Beltrami
operator on $\Gamma$, in this order. We then consider the initial-boundary value problem
\beqa
  && \bigl(\coeff \bigr) \, \pt\mu
  + \mu \, g'(\rho) \, \pt\rho
  - \Delta\mu =0\quad\mbox{in $Q$,}
  \label{ss1}
  \\
	&&\pn\mu =0 \quad\mbox{on $\Sigma$},\label{ss2}\\
  && \pt\rho - \Delta\rho + \xi+\pi(\rho) = \mu \, g'(\rho)	
	\quad\mbox{in $Q$},
  \label{ss3}
  \\
	&&\xi\in\beta(\rho)\quad\mbox{a.\,e. in $Q$},\label{ss4}\\
  && \pn\rho+\pt\rg+\xig+
	\pig(\rg) -\Dg	\rg= \ug,
	\quad \rg =\rho_{\pier{|\Sigma}},\quad \mbox{on $\Sigma$},
  \label{ss5}
  \\
	&&\xig\in\beg(\rg)\quad\mbox{a.\,e. on $\Sigma$},\label{ss6}\\
  && \mu(0) = \mu_0,\quad \rho(0)=\rho_0,
	\quad\mbox{in $\Omega$}, \quad
    \gianni\rg(0) = \rho_{0_{|\Gamma}}\quad\mbox{on $\Gamma$.}
  \label{ss7}
\eeqa
\pier{We point out that $\beta$ and $\beta_\Gamma$ denote two maximal monotone graphs, with domains $D(\beta) \supseteq D(\beta_\Gamma)$ that are intervals containing $0$; $\beta$ and $\beta_\Gamma$ fulfill suitable properties (cf. the later {\bf (A3)} and {\bf (A7)}) about their values and growths; $\pi $ and $\pig$ denote two Lipschitz continuous perturbations;
$g$ is a smooth nonnegative and concave function defined on $D(\beta) $; $\ug$ is a boundary datum.}
 
\pier{The system \eqref{ss1}--\eqref{ss7} is related to a model for phase segregation through atom rearrangement on a lattice that has been proposed in~\cite{PG}. This model (see also \cite{CGPS3} for a detailed derivation) is a modification of the Fried--Gurtin approach to
phase segregation processes (cf. \cite{FG,G}). The order parameter 
$\rho$, which in many cases represents the (normalized) density of one of the phases, 
and the chemical potential $\mu$ are the unknowns of the system.}

\pier{Let us note at once that, in our present contribution, the Neumann homogeneous boundary condition \eqref{ss2} is considered for $\mu$, while, differently from previous analyses\cite{CGKPS, CGPS3, CGPS7,CGPSco,CGPS6,CGPS4,CGPS5,CGSco1}, the dynamic boundary condition \eqref{ss5} is assumed for $\rg$, the trace of $\rho$ on the boundary.}

\pier{The approach by Podio-Guidugli~\cite{PG} is based on a local 
free energy density (in the bulk) of the form
\begin{equation}
\label{fe1}
\psi(\rho,\nabla\rho,\mu)=-\mu\,\rho+W(\rho)+\frac 1 2\,|\nabla\rho|^2,
\end{equation}
where $W$ is a double-well potential, whose derivative (in the differentiable case) plays as the sum $\beta + \pi $ in \eqref{ss3}
(see also \eqref{ss4}). By \eqref{fe1}, one arrives at the evolutionary system
\begin{align}
  & 2\rho \, \partial_t\mu
  + \mu \, \partial_t\rho
  - \Delta\mu = 0
  \label{oldprima}
  \\
  & - \Delta\rho + W'(\rho) = \mu\,. 
   \label{oldseconda}
\end{align}
The above equations are assumed to hold in $Q$ and must be complemented 
with boundary and initial conditions.
The typical example for} \juerg{a} \pier{smooth double-well potential~$W$ is
\begin{equation}
  W (r)= \frac14(r^2-1)^2 \,,
  \quad r \in \rz ,
  \label{regpot}
\end{equation} 
while another smooth potential, but defined on a bounded interval, is given by  
\begin{equation} 
 W(r) = (1+r)\ln (1+r)+(1-r)\ln (1-r) - c r^2 \,,
  \quad r \in (-1,1),
  \label{logpot}
\end{equation}
where the coefficient $c$ is taken greater than $1$ in order that $W$ be nonconvex. 
The potentials \eqref{regpot} and \eqref{logpot}
are usually referred to as the {\em classical regular\/}
and the {\em logarithmic double-well\/} potential, respectively.
Observe that the derivative of the logarithmic potential
becomes singular at $\,\pm 1$.
However, one can consider nondifferentiable potentials,
where an important example is given by the so-called 
{\em double-obstacle} potential
\beq
  W(r) = I_{[-1,1]}(r) - c r^2 \,,
  \quad r \in \rz,
  \label{obspot}
\eeq
where $c>0$ is a positive constant and $I_{[-1,1]}:\rz\to[0,+\infty]$ denotes the indicator function of~$[-1,1]$, i.e.,
we have} \juerg{$I_{[-1,1]}(r)=0$ if $|r|\leq1$ and $I_{[-1,1]}(r)=+\infty$} \pier{otherwise.
In this case, the order parameter is subjected to the unilateral constraint $|\rho|\leq1$,
and \eqref{oldseconda} should be read as a differential inclusion where $W'(\rho) = \beta(\rho) + \pi (\rho)$, with $\beta$  representing the
subdifferential $\partial I_{[-1,1]}$ of $I_{[-1,1]}$ and $\pi (\rho)= -2c\rho.$}

\pier{The system \eqref{oldprima}--\eqref{oldseconda} is a variation of the
Cahn--Hilliard system originally introduced in \cite{CahH}
and first studied mathematically in the seminal paper \cite{EllSh} (for an updated
list of references on the Cahn--Hilliard system, see~\cite{Heida}). 
An initial-boundary value problem for \eqref{oldprima}--\eqref{oldseconda} is in general ill-posed: indeed, as it was pointed out
in \cite{CGPS6}, when assuming Neumann homogeneous boundary
conditions for both $\rho$ and $\mu$, the related problem  
may have infinitely many smooth and even nonsmooth solutions. 
Then, two small regularizing parameters $\varepsilon>0$ and $\delta>0$ were considered in 
\cite{CGPS3}, which led to the regularized model equations 
\begin{align}
  & \bigl( \eps + 2\rho \bigr) \, \partial_t\mu
  + \mu \, \partial_t\rho
  - \,\Delta\mu = 0,
  \label{viscprima}
  \\
  & \delta\, \partial_t\rho - \Delta\rho + F'(\rho) = \mu\,. 
   \label{viscseconda}
\end{align}
The system \eqref{viscprima}--\eqref{viscseconda} constitutes
a modification of the so-called {\em viscous} Cahn--Hilliard system (see\cite{NC} and the recent papers\cite{CGW, CGS0, CGS2} along with their references). We point out that
 \eqref{viscprima}--\eqref{viscseconda} was analyzed, in the case of 
no-flux boundary conditions for both $\mu$ and $\rho$, in the 
papers \cite{CGPS3,  CGPSco, CGPS4, CGSco1} concerning
well-posedness, regularity, asymptotics as $\eps \searrow 0$, and optimal control. Later, the
local free energy density \eqref{fe1} was generalized to the form
\begin{equation}
\label{fe2}
\psi(\rho,\nabla\rho,\mu)=-\mu\,g(\rho)+W(\rho)+\frac 1 2\,|\nabla\rho|^2,
\end{equation}
{with a function $g$ having suitable properties. If one puts,
without loss of generality, $\varepsilon=\delta=1$, then one obtains the
more general system}
\begin{align}
  & {\bigl( \coeff \bigr) \, \partial_t\mu
  + \mu \, g'(\rho) \, \partial_t\rho
  - \Delta\mu = 0,}
  \label{locprima}
  \\
  & {\partial_t\rho - \Delta\rho + W'(\rho) = \mu \,g'(\rho)},
   \label{locseconda}
\end{align}
which was investigated in the contributions 
\cite{CGKPS, CGPS7,CGPS6,CGPS5}, still for no-flux boundary conditions, 
also from the side of the numerical approximation. The related phase relaxation system
(in which the diffusive term $- \Delta\rho $ disappears from \eqref{locseconda}), has been dealt with in \cite{CGKS1, CGKS2, CGS-lom}. We also mention the recent article \cite{CGS3},
where a nonlocal version of \eqref{locprima}--\eqref{locseconda} -- based on the replacement of the diffusive term of \eqref{locseconda} with a nonlocal operator acting on $\rho$ -- has been largely investigated.}

\pier{In the present paper, we consider a total free energy of the system which also includes 
a contribution on the boundary; in fact, we postulate that a phase transition phenomenon is occurring as well on the boundary, and the physical variable on the boundary is just the trace of the phase variable in the bulk. Then, we choose
a total free energy functional of the form
\begin{align}
\label{fe3}
\mathbf{\Psi} [\mu(t), \rho\juerg{(t)}, \rg\juerg{(t)} ]
= &\int_\Omega \Bigl[-\mu\,g(\rho)+W(\rho)+\frac 1 2\,|\nabla\rho|^2\Bigr]
\juerg{(t)\dx} \non\\
&+  \int_\Gamma \Bigl[-u_\Gamma\,\rg +W_\Gamma (\rho_\Gamma)+\frac 1 2\,|\nabla_\Gamma \rho_\Gamma|^2\Bigr]\juerg{(t)}\,\dg , \juerg{\quad t\in [0,T],}
\end{align}
where $W_\Gamma$ denotes a double-well potential \juerg{having more or less the same 
behavior as} $W$, 
and $u_\Gamma$ is a source term \juerg{that} may exert a (boundary) control on the system. From this expression of the total free energy, one recovers the PDE system resulting in 
\eqref{ss1}--\eqref{ss7}; in particular, we point out that the derivative or subdifferential of $W_\Gamma$ is expressed by the sum $\beg + \pig$ in \eqref{ss5}--\eqref{ss6}.}

\pier{Our aim here is \juerg{to prove} the well-posedness of \eqref{ss1}--\eqref{ss7} as well as regularity properties of the solution, like the $L^\infty$-boundedness of both variables $\mu$ and $\rho$,
and consequently of $\rg$ on the boundary, and the so-called strict separation property, in the case when $W$ and $W_\Gamma$ behave like \eqref{logpot}, \juerg{which means to show}  that $\rho$
and $\rg$ are uniformly bounded \juerg{away} from $-1$ and $+1$. We employ certain techniques which combine some key ideas essentially from the papers\cite{CGPS3} and \cite{CGPS5} together with the treatment 
of dynamic boundary conditions devised in \cite{Calcol} and exploited in other 
solvability studies and optimal control theories, namely, 
\cite{CF1,CF2,CF3,CGS0,CGS1,CGS2,CS}. For these reasons, we often refer to 
the abovementioned papers in running the proofs.} 

\pier{About dynamic boundary conditions, let us notice that there has been a recent growing interest about the justification and the study of phase field models, as well as systems of Allen--Cahn and Cahn--Hilliard type, 
\emph{including dynamic boundary conditions}. Without trying to be exhaustive, we 
\juerg{mention} at least the papers\cite{ChGaMi, CGM1, CGM2, GaGu, GaWa, GiMiSc, GMS10,  GoMi, GMS11, Is, Li, MRSS}.}

The paper is organized as follows: in Section 2, we formulate the relevant assumptions on the
data of system \eqref{ss1}--\eqref{ss7}, and we state the main results of this paper. Section 3 
then brings a detailed proof of the existence result stated in Theorem 2.1, while Section 4
deals with the proofs of Theorem 2.2 and Theorem 2.4. 

Throughout the paper, \pier{for a general Banach space $X$, we denote} by $\|\cdot\|_X$ its norm and by $X'$ its dual space.
The only exemption from this convention are the norms of the $L^p$ spaces and of  their powers, which we often denote by
$\|\cdot\|_p$, for $1\le p\le +\infty$.   
Moreover, we often utilize the continuity of the embedding $\Huno \subset\Lx p$
for $1\leq p\leq 6$ and the related Sobolev inequality
\beq
  \|v\|_p \leq \gianni{C_\Omega} \|v\|_{\Huno}
  \quad \hbox{for every $v\in \Huno$ and $1\leq p \leq 6$,}
  \label{sobolev}
\eeq
where $\gianni{C_\Omega}$ depends only on~$\Omega$. Notice that these embeddings  
are compact for $1\le p<6$. We also use the corresponding compactness inequality
\beq
  \|v\|_4 \,\leq\, \delta\,\|\nabla v\|_2 + \widetilde C_\delta \|v\|_2
  \quad \hbox{for every $v\in \Huno$ and $\delta>0$,}
  \label{compact}
\eeq
where $\widetilde C_\delta$ depends only on $\Omega$ and~$\delta$, and recall that
the embedding \,$\Hdue\subset C^0(\overline{\oma})$ is compact.
Furthermore, we make repeated use of the notation
\beq
  Q_t := \Omega \times (0,t),\quad\Sigma_t:=\Gamma\times (0,t),
  \quad \hbox{for $t\in (0,T],$}
  \label{defQt}
\eeq
as well as of H\"older's inequality and of the elementary Young inequalities
\begin{align}
\label{Young}
& |ab| \leq \gamma\, |a|^2 + \frac 1{4\gamma} \, |b|^2\,\,\quad\mbox{and}\quad\,\,
  |ab|\leq \frac {\gamma^p}p\,|a|^p+\frac{\gamma^{-q}}q\,|b|^q,\non\\[1mm]
&\hbox{for every $a,b\in \rz$, $\gamma>0$, and $1<p,q<+\infty$ with $\frac 1p+\frac 1q=1$}.
\end{align}

\section{General assumptions and main results}
\setcounter{equation}{0}

In this section, we formulate the general assumptions for the data of the system 
\eqref{ss1}--\eqref{ss7}, and we state the main results of this paper. To begin with, we
introduce some denotations. We set
\begin{align*}
&H:=\Ldue,\quad V:=\Huno, \quad W:=\{w\in\Hdue:\mbox{\,$\pn w=0$\, on $\,\Gamma$}\},\\
&\Hg:=L^2(\Gamma),\quad\Vg:=H^1(\Gamma),
\end{align*}
and endow these spaces with their standard norms. Notice that we have $V\subset H\subset V'$ and 
$\Vg\subset\Hg\subset \Vg'$ with dense, continuous and compact embeddings. Moreover, for every
proper, convex and lower semicontinuous function $\,\pier{\widehat\alpha{}}:\rz\to [0,+\infty]\,$ satisfying $\,\pier{\widehat\alpha{}}(0)=0$,
we denote by $\,\,\alpha:=\partial\pier{\widehat\alpha{}}\,\,$ its subdifferential, which is known to be a 
maximal monotone graph in $\rz\times\rz$ satisfying $\,0\in\alpha(0)$. We denote its effective domain by
$\,D(\alpha)\,$, and, for $\,r\in D(\alpha)$, by $\,\alpha^\circ(r)$\, the element
of $\,\alpha(r)\,$ having minimal modulus. Moreover, for any $\eps>0$, we denote
by \,$\alpha^\eps_Y$\, the Yosida approximation of $\,\alpha\,$ 
at the level $\eps>0$. As is well known (see, e.\,g, \cite[p.~28]{Brezis}), $\alpha^\eps_Y\,$ is 
a monotone and Lipschitz continuous function on $\rz$, and we have that
\beq\label{domeps}
|\alpha^\eps_Y(r)|\le|\alpha^\circ(r)|\quad\forall\,\eps>0,\,\quad\mbox{as well as}\,\quad
\lim_{\eps\searrow0}\,\alpha_Y^\eps(r)=\alpha^\circ(r),\quad\,\,\mbox{for all \,$r\in D(\beta)$}. 
\eeq
Moreover, the antiderivative
\beq\label{antider1}
\pier{\widehat\alpha{}}^\eps(r):=\int_0^r\alpha^\eps_Y(s)\ds, \quad r\in\rz,
\eeq
is a convex function on $\rz$ that satisfies 
\beq\label{antider2}
0\le \pier{\widehat\alpha{}}^\eps_Y(r)\le\pier{\widehat\alpha{}}(r),\quad\pier{\widehat\alpha{}}^\eps_Y(r)\nearrow \pier{\widehat\alpha{}}(r) 
\mbox{\, as $\,\eps\searrow 0\,$, \,\, for all \,$r\in\rz$}.
\eeq

We make the following general assumptions:

\vspace{3mm}\noindent
{\bf (A1)} \quad $\mu_0\in W$, $\,\,\,\mu_0\ge 0$ a.\,e. in $\oma$,
\,\,$\,\rho_0\in \Hdue$, \,\,\,$\rho_{0_\Gamma}:=\rho_{0_{|\Gamma}}\in H^2(\Gamma)$.

\vspace{1mm}\noindent
{\bf (A2)} \quad $\ug\in H^1(0,T;\Hg)$.

\vspace{1mm}\noindent
{\bf (A3)} \quad  $\beta=\partial\pier{\widehat\beta{}},\,\,\,\beta_\Gamma=\partial\pier{\widehat\beta{}}_\Gamma$, \,where
\,\,$\pier{\widehat\beta{}},\pier{\widehat\beta{}}_\Gamma:\rz\to [0,+\infty]$\, are proper, convex, and \,lower 
\\ \hspace*{15.5mm}semicontinuous
functions satisfying $\pier{\widehat\beta{}}(0)=\pier{\widehat\beta{}}_\Gamma(0)=0$. 

\vspace{2mm}\noindent
We remark that {\bf (A3)} entails that $\beta$ and $\beta_\Gamma$ are 
maximal monotone graphs in $\rz\times\rz$, with $0\in\beta(0)$ and $0\in \beta_\Gamma(0)$, whose
Yosida approximations $\beta^\eps_Y$ and $\beta^\eps_{\Gamma_Y}$ satisfy the conditions
\eqref{domeps}--\eqref{antider2}. In particular, we have $\beta^\eps_Y(0)=\beta^\eps_{\Gamma_Y}(0)=0$. 

\vspace{2mm}\noindent 
{\bf (A4)} \quad The functions $\pi,\pig:\rz\to\rz$ are Lipschitz continuous.

\vspace{1mm}\noindent 
{\bf (A5)} \quad The function $\,g:\overline{D(\beta)}\to [0,+\infty)$ belongs to
$C^2$, is bounded and concave, and\\ \hspace*{15.5mm}$g'$ is bounded and Lipschitz continuous.

\vspace{1mm}\noindent
{\bf (A6)} \quad $\pier{\widehat\beta{}}(\rho_0)\in L^1(\oma)$, \,\,\,$\pier{\widehat\beta{}}_\Gamma(\rho_{0_\Gamma})
\in L^1(\Gamma)$, \,\,\,$\beta^\circ(\rho_0)\in H$, \,\,\,$\beta_\Gamma^\circ(\rho_{0_\Gamma})
\in \Hg$.
 
\vspace{3mm}\noindent Finally, we need a compatibility condition, which essentially means that
the graph in the bulk is dominated by the graph on the boundary:

\vspace{1mm}\noindent
{\bf (A7)} \quad $D(\beta_\Gamma)\subset D(\beta)$, and there are constants $\,\eta>0\,$ and
$\,C_\Gamma\ge 0\,$ such that
\beq
\label{domination1}
|\beta^\circ(r)|\,\le\,\eta\,|\beta_\Gamma^\circ(r)|\,+\,C_\Gamma \quad\forall\,
r\in D(\beta_\Gamma).
\eeq  

We now state the main results of this paper. Concerning \gianni{well-posedness}, we have the following result.

\vspace{5mm}\noindent
{\sc Theorem 2.1:}  \quad\,\,{\em Suppose that the assumptions} {\bf (A1)}--{\bf (A7)} {\em are 
fulfilled. Then the system} \eqref{ss1}--\eqref{ss7} {\em admits a unique solution $\,(\mu,\rho,\rg,\xi,\xig)\,$
such that \,$\mu\ge 0\,$ almost everywhere in $\,Q$\, and}
\begin{align}
&\mu\in W^{1,p}(0,T;H)\cap C^0([0,T];V)\cap L^p(0,T;W)\cap L^\infty(Q) \quad\mbox{for every 
\,$p\in [1,+\infty)$},
\label{regmu}\\ 
&\rho\in W^{1,\infty}(0,T;H)\cap H^1(0,T;V)\cap L^\infty(0,T;\Hdue),\label{regrho}\\
&\rg\in W^{1,\infty}(0,T,\Hg)\cap H^1(0,T,\Vg)\cap L^\infty(0,T;H^2(\Gamma)),\label{regrg}\\
&\xi\in L^\infty(0,T;H),\quad\xig\in L^\infty(0,T;\Hg).\label{regxi}
\end{align}

\vspace{5mm}
\noindent
{\sc Remark 2.2:} \quad\,\,Observe that it follows from standard embedding results (see, e.\,g., 
\cite[Sect.~8, Cor.~4]{Simon}) that $\rho\in C^0([0,T];H^s(\oma))$ for $0<s<2$. In particular,
$\,\rho\in C^0(\overline{Q})$, which entails that $\,\rg=\rho_{\pier{|\Sigma}}\in C^0(\overline{\Sigma})$.
Notice also that it follows from \eqref{regmu} and \eqref{regrho} that both $\,|\nabla\mu|\,$
and $\,|\nabla\rho|\,$ belong to $L^4(0,T;L^6(\oma))$.

\vspace{3mm}
Under somewhat stronger assumptions, we can \gianni{prove that the solution enjoys a better regularity}. 
More precisely, we have the following result.

\vspace{5mm}\noindent
{\sc Theorem 2.3:} \quad\,\,{\em Suppose that} {\bf (A1)}--{\bf (A7)} {\em hold true, and let
\beq
  \gianni{%
  \ug\in L^\infty(\Sigma), \quad
  \beta^\circ(\rho_0)\in L^\infty(\oma)
  \aand
  \beta^\circ(\rho_{0_\Gamma})\in L^\infty(\Gamma).}%
  \label{hpxibdd}
\eeq
Then, the unique solution $\,(\mu,\rho,\rg,\xi,\xig)\,$ to \eqref{ss1}--\eqref{ss7}
established in Theorem~2.1 also satisfies  $\,\xi\in L^\infty(Q)$.}

\vspace{3mm}

\gianni{The above theorems hold for every pair of potentials 
\juerg{$\widehat\beta$} and $\pier{\widehat\beta{}}_\Gamma$
satisfying assumptions {\bf (A3)} and~{\bf (A7)}.
In our last result, we reinforce the compatibility condition a little
and consider a more restricted class of potentials
which does not include potentials of obstacle type,
while potentials of logarithmic type are still admitted. 
Namely, we assume~that}

\vspace{3mm}\noindent
\gianni{{\bf (A8)} \quad 
$-\infty\leq r_-<r_+\leq+\infty$\aand $D(\beta)=D(\beta_\Gamma)=(r_-,r_+)$}

\vspace{3mm}\noindent
\gianni{and prove a uniform separation property.
We observe that such an assumption allows multi-valued subdifferentials.
A~stability estimate holds if in addition the \pier{potentials} are smooth.}

\vspace{5mm}\noindent
\gianni{{\sc Theorem 2.4:} \quad\,\,{\em Assume that} {\bf (A1)}--{\bf (A8)} 
{\em and} \eqref{hpxibdd} {\em are fulfilled}.
{\em Then, the unique solution $\,(\mu,\rho,\rg,\xi,\xig)\,$ to \eqref{ss1}--\eqref{ss7}
established in Theorem~2.1 satisfies the 
uniform separation property
\beq\label{separation}
r_*\le\rho(x,t)\le r^*,\quad \mbox{\emph{for\,\,\,every \,\,\,}}(x,t)\in \overline{Q},
\eeq
with constants $\,r_*,r^*\in(r_-,r_+)\,$ that depend only on the data of the system,
and we also have $\xi_\Gamma\in L^\infty(\Sigma)$.
Moreover, assume that $\pier{\widehat\beta{}}$ and $\pier{\widehat\beta{}}_\Gamma$ 
are functions of class $C^2$ in $(r_-,r_+)$
and let $\,u_{\pier{i_\Gamma}}\in H^1(0,T;\Hg)\cap L^\infty(\Sigma)\,$ be given.
Then, the corresponding solutions 
$(\mu_i, \rho_i,\rho_{i_\Gamma},\xi_i,\xi_{i_\Gamma})$
satisfy the stability estimate}}
\begin{align}
\label{stabu}
&\|\mu_1-\mu_2\|_{L^\infty(0,t;H)\cap L^2([0,t];V)}\,+\,
\|\rho_1-\rho_2\|_{H^1(0,t;H)\cap C^0([0,t];V)\cap L^2(0,t;\Hdue)}\non\\
&\ +\|\rho_{1_\Gamma}-\rho_{2_\Gamma}\|_{H^1(0,t;\Hg)
\cap C^0([0,t];\Vg)\cap L^2(0,t;H^2(\Gamma))}\non\\
&{}\le\,\widehat C\,\|u_{\pier{1_\Gamma}}-u_{\pier{2_\Gamma}}\|
_{L^2(0,t;\Hg)}\,,\quad\forall\,t\in (0,T],
\end{align}
{\em with a finite constant $\widehat C>0$ \gianni{that depends  only on the data of the system}.}

\vspace{5mm}
\noindent
{\sc Remark 2.5:} \quad
\gianni{If $(r_-,r_+)$ is bounded
(but a similar remark holds in the case of a half line),
then both $\beta$ and $\beta_\Gamma$ become singular at~$r_\pm$
due to maximal monotonicity.
However, such a singularity never becomes active thanks to~\eqref{separation}.}


\section{Existence \gianni{and uniqueness}}
\setcounter{equation}{0}
In this section, we prove the relevant \gianni{well-posedness} result.

\vspace{2mm}\noindent
{\sc Proof of Theorem 2.1:}\quad We employ  an 
approximation scheme based on Yosida approximation and a time delay in the right-hand side
of \eqref{ss3}. To this end, let
$\,\beta^\eps:=\beta_Y^\eps\,$ and $\,\beta_\Gamma^\eps:=\beta_{\Gamma_Y}^{\eps\eta}$ denote the Yosida
approximations of $\beta$ and $\beta_\Gamma$ at the levels $\eps>0$ and $\eps\eta>0$, respectively,
where $\eta>0$ is the constant introduced in \eqref{domination1}. Then $\beta^\eps\,,\,\beta_\Gamma^\eps\,$
satisfy $\,\beta^\eps(0)=\beta_\Gamma^\eps(0)=0$ and the conditions \eqref{domeps}--\eqref{antider2}
correspondingly. Moreover, we infer from \cite[Lemma~4.4]{Calcol} that
\beq\label{domination2}
|\beta^\eps(r)|\,\le\,\eta\,|\beta_\Gamma^\eps(r)|\,+\,C_\Gamma \quad
\forall\,r\in\rz.
\eeq 
Since, by virtue of their monotonicity, the functions $\beta^\eps$ and $\beta^\eps_\Gamma$ have equal sign,
it follows from \eqref{domination2} and Young's inequality for all $r\in\rz$ the estimate
\begin{align}
\label{domination3}
&\beta^\eps(r)\,\beta_\Gamma^\eps(r)\,=\,|\beta^\eps(r)|\,|\beta_\Gamma^\eps(r)|
\,\ge\,\frac 1\eta\,|\beta^\eps(r)|^2-\frac{C_\Gamma}\eta\,|\beta^\eps(r)| 
\,\ge\,\frac 1{2\eta}\,|\beta^\eps(r)|^2-\frac{C_\Gamma^2}{2\eta}\,.
\end{align}
We also temporarily extend the function $g$ to the whole real line $\rz$, still terming
the extended function $g$, in such a way that
\begin{align}\label{extendg}
&\mbox{$g \in C^1(\rz)$, $g$ and $g'$ are bounded and Lipschitz continuous on $\rz$, and}\non\\
&g(r)\ge -1/3 \quad\mbox{(i.\,e., $1+2g(r)\ge 1/3>0$)} \quad\forall\,r\in\rz.
\end{align}

To \gianni{introduce the time delay}, we define for every $\tau\in (0,T)$ the translation operator 
$\,\T:C^0([0,T];H)\to C^0([0,T];H)$ by setting, for all $v\in C^0([0,T];H)$,
\beq
\label{defT}
\T(v)(t):=v(t-\tau) \quad\mbox{if $\,t>\tau$}\quad\mbox{and }\,\T(v)(t):=v(0) \quad
\mbox{if }\,t\le\tau.
\eeq
Notice that for every $v\in H^1(0,T;H)$ it holds that
\begin{align}
\label{T1}
\|\T(v)\|_{L^2(Q_t)}^2\,&\le\,\left\{
\begin{array}{l}
\|v\|_{L^2(Q_t)}^2\,+\,\tau\,\|v(0)\|^2_H
\quad\mbox{for all $t\in [\tau,T]$},\\[1mm]
t\,\|v(0)\|_H^2\quad\mbox{for all $t\in [0,\tau]$},
\end{array}
\right. \\[2mm]
\label{T2}
\|\pt\T(v)\|_{L^2(Q_t)}^2\,&\le\,\|\gianni{\pt v}\|_{L^2(Q_t)}^2\quad\mbox{for a.\,e. 
$\,t\in (0,T)$},
\end{align}
while for every $v\in C^0([0,T];V)$ we have 
\begin{align}
\label{T3}
\|\nabla\T(v)\|_{L^2(Q_t)}^2\,&\le\,\left\{
\begin{array}{l}
\|\nabla v\|_{L^2(Q_t)}^2\,+\,\tau\,\|\nabla v(0)\|^2_H
\quad\mbox{for all $t\in [\tau,T]$},\\[1mm]
t\,\|\nabla v(0)\|_H^2\quad\mbox{for all $t\in [0,\tau]$}.
\end{array}
\right. 
\end{align}

We then consider the problem (which in the following will be termed $(P_{\tau}^\eps)$) of 
finding a triple $(\mute,\rhote,\rhotge)$ with
\begin{align}
\label{regepstau}
&\mute\in H^1(0,T;H)\cap C^0([0,T];V)\cap L^2(0,T;W)\cap L^\infty(Q),
\quad \mute\ge 0\quad\aeQ,\non\\
&\rhote\in W^{1,\infty}(0,T;H)\cap H^1(0,T;V)\cap L^2(0,T;\Hdue),\non\\
&\rhotge\in W^{1,\infty}(0,T;H_\Gamma)\cap H^1(0,T;V_\Gamma))\cap L^2(0,T;
H^2(\Gamma)),
\end{align}
which solves the initial-boundary value problem
\begin{align}
  & \bigl(1 + 2 g(\rhote)\bigr)\, \pt\mute + g'(\rhote)\,\pt\rhote\,\mute- \Delta\mute = 0
   \quad \mbox{a.\,e. in $Q$},
  \label{etau1}
  \\
   \label{etau2}
  &\pn\mute=0\quad\mbox{a.\,e. on $\Sigma$},\\
  \label{etau3}
  & \pt\rhote - \Delta\rhote + \beta^\eps(\rhote)+\pi(\rhote) = {\cal T}_\tau(\mute)\,g'(\rhote)
  \quad \mbox{a.\,e. in $Q$},\\
  \label{etau4}
  &\pn\rhote+\pt\rhotge+\beta^\eps_\Gamma(\rhotge)+
	\pig(\rhotge) -\Dg	\rhotge= \ug,
	\quad \rhotge=\rho^\eps_{\tau_{|\Sigma}},\quad \mbox{a.\,e. on $\Sigma$},\\
  \label{etau5}
   & \mute(0) = \mu_0
  \quad\mbox{and}\quad
  \rhote(0) = \rho_0
  \quad \mbox{a.\,e. in $\oma$}, \quad \rhotge(0)=\rho_{0_{|\Gamma}}\quad\mbox{a.\,e. on $\Gamma$}.
  \end{align}
For convenience, we allow $\tau$ to attain just discrete values,
namely, $\tau=\tau_N:=\frac T N$, where $N$ is any positive integer.
We now proceed in four steps: at first, we show that the problem $(P_\tau^\eps)$
is well posed for every $\tau=\tau_N$, $N\in\nz$, then we prove a number of a priori
estimates, followed by the limit processes as $N\to\infty$ and then as $\eps\searrow0$.
The limit processes use standard compactness and monotonicity arguments.
As a notational convention, in the remainder of this proof \pier{we will}
denote by $\,K_i$, $i\in\nz$, positive constants that may depend on the
data of the system \eqref{ss1}--\eqref{ss7} but neither on $\tau>0$ nor on $\eps>0$.

\vspace{5mm}\noindent
\underline{\sc Step 1:}\quad We claim that, for every $\tau_N=\frac T N$, $N\in\nz$, and every $\eps>0$, the problem \eqref{etau1}--\eqref{etau5}
admits a unique solution triple $\bigl(\mu_{\tau_N}^\eps,\rho_{\tau_N}^\eps,
\rho^\eps_{\tau_{N_\Gamma}}\bigr)$ satisfying \eqref{regepstau}.

To prove this claim, we set $t_n:=n\tau_N$ for $n=0,\dots,N$, and observe
that the problem of solving \eqref{etau1}--\eqref{etau5} 
becomes equivalent to a finite sequence of $N$ problems
that can be solved inductively step by step.
However, instead of considering the natural time intervals
$[t_{n-1},t_n]$, $n=1,\dots,N$ and glueing the solutions 
together, we solve $N$ problems
on the time intervals $I_n=[0,t_n]$, $n=1,\dots,N$,
by constructing the solution directly on the whole of~$I_n$ at each step.
These problems read as follows:
\begin{align}
  & (1 + 2g(\rhon)) \,\pt\mun + g'(\rhon)\,\pt\rhon\,\mun- \Delta\mun = 0
  \aand
  \mun \geq 0
  \quad \hbox{a.\,e.\ in $\Omega\times I_n$},
  \label{sn1}
  \\
  & \pn\mun = 0 \quad \hbox{a.\,e. on $\Gamma\times I_n$}
  \aand
  \mun(0) = \mu_0 \quad\mbox{a.\,e. in $\oma$},
  \label{sn2}
  \\
  \label{sn3}
  &\pt\rhon - \Delta\rhon + \beta^\eps(\rhon)+\pi(\rhon) = {\cal T}_{\tau_N}(\munmu)\,g'(\rhon)
  \quad\hbox{a.\,e.\ in $\Omega\times I_n$},\\[0.1cm]
  \label{sn4}
  &  \pn\rhon+\pt\rhong+\beta_\Gamma^\eps(\rhong)+
	\pig(\rhong) -\Dg	\rhong= \ug, \non \\
	& \hskip5cm \rhong=\rho_{n_{|\Gamma\times I_n}},\quad \mbox{a.\,e. on $\Gamma\times I_n$},\\[0.1cm]
	\label{sn5}
   & \rhon(0) = \rho_0\quad\mbox{a.\,e. in }\,\oma, \quad \rhong(0)=\rho_{0_{|\Gamma}}
   \quad\mbox{a.\,e. on }\,\Gamma.
\end{align}
Their (unique) solutions $(\mun,\rhon,\rhong)$ are required to satisfy the regularity properties inductively
obtained by taking $t_n$ in place of $\,T\,$ in \eqref{regepstau}. 
The operator ${\cal T}_{\tau_N}$ appearing on the \rhs\ of~\eqref{sn3} 
acts on functions that are not defined in the whole of~$(0,T)$; however,
its meaning is still given by~\eqref{defT} if $n>1$,
while for $n=1$ we simply set ${\cal T}_{\tau_N}(\munmu)=\mu_0$.

Observe that the unique solvability of \eqref{sn1}--\eqref{sn5} for every 
$n\in\{1,\ldots,N\}$ would imply that 
$$
\rho_{N_{|\oma\times I_{N-1}}}=\rho_{N-1},\quad \mu_{N_{|\oma\times I_{N-1}}}=\mu_{N-1},
$$
which would entail that 
\begin{align*}
{\cal T}_{\tau_N}(\mu_{N-1})(x,t)=\mu_{N-1}(x,t-\tau_N)=\mu_N(x,t-\tau_N)=
{\cal T}_{\tau_N}(\mu_N)(x,t) \\
 \quad\mbox{for a.\,e. $(x,t)\in \pier{\oma\times I_{N-1}}$}.
\end{align*}
It then would follow that $(\mu_N,\rho_N,\rho_{N_\Gamma})$ is the unique solution 
to the system \eqref{etau1}--\eqref{etau5} for $\tau=\tau_N$. 
Consequently, it suffices to show that the system  \eqref{sn1}--\eqref{sn5} 
enjoys for every $n\in\{1,\ldots, N\}$ 
a unique solution having the required
regularity. To this end, we argue by induction. Since the proof for $n=1$ is similar to 
that used in the induction step $n-1\longrightarrow n$ (recall that we have set
${\cal T}_{\tau_N}(\mu_0)=\mu_0$ on $\oma\times (0,t_1)$), we may confine 
ourselves to just
perform the latter.

So let $1<n\le N$, and assume that, for $1\le k\le n-1$, unique solutions 
$(\mu_k,\rho_k,\rho_{k_\Gamma})$ to the system \eqref{sn1}--\eqref{sn5} have already been 
constructed that
satisfy for $1\le k\le n-1$ the conditions
\begin{align}
\label{regtau}
&\mu_k\in H^1(I_k;H)\cap C^0(\bar I_k;V)\cap L^2(I_k;W)
\cap L^\infty(\oma\times I_k),\quad \mu_k\ge 0\quad
\mbox{a.\,e. in $\oma\times I_k$},\non\\
&\rho_k\in W^{1,\infty}(I_k;H)\cap H^1(I_k;V)\cap L^2(I_k;\Hdue),\non\\
&\rho_{k_\Gamma}\in W^{1,\infty}(I_k;H_\Gamma)\cap H^1(I_k;V_\Gamma))\cap L^2(I_k;
H^2(\Gamma)).
\end{align}

Now observe that the system \eqref{sn3}--\eqref{sn5} is nothing but the strong
formulation of the variational problem \cite[Eqs.~(4.5)--(4.8)]{Calcol} (if one puts, in the
notation used there, $f_\Gamma:=u_\Gamma, f=0$), with the only exception that the
expression $\pi(\rhone)$ occurring there is in our case replaced by $\,\pi(\rhone)
-{\cal T}_{\tau_N}(\munmu)\,g'(\rhone)$. 
Now notice that the data of the system satisfy all of the conditions imposed in \cite{Calcol}.
In view of \eqref{extendg}, and since
${\cal T}_{\tau_N}(\munmu)$ obviously belongs to $L^\infty(\oma\times I_n)$,
we can therefore adopt the arguments employed in \cite[Sect.~4]{Calcol} (see, in particular, 
\cite[Prop.~4.1 and Rem.~4.5]{Calcol}) to conclude that \eqref{sn3}--\eqref{sn5}
enjoys a unique solution $(\rhon,\rhong)$ such that
\begin{align}
\label{est1}
&\rhon\in H^1(I_n;H)\cap C^0(\bar I_n;V)\cap L^2(I_n;\Hdue),\\
&\rhong\in H^1(I_n;H_\Gamma)\cap C^0(\bar I_n;V_\Gamma)\cap L^2(I_n;
H^2(\Gamma))\\
&\beta^\eps(\rhon)\in L^2(I_n;H),\quad \beta^\eps_\Gamma(\rhong)\in
L^2(I_n;H_\Gamma)\,.
\end{align}
In order to recover the full regularity required in \eqref{regepstau},
we still need to show that 
\beq\label{extrareg}
\pt\rhon\in L^\infty(I_n;H)\cap L^2(I_n;V), \quad \pt\rhong \in L^\infty(I_n;\Hg)
\cap L^2(I_n;\Vg).
\eeq
To establish \eqref{extrareg}, we proceed formally, noting that the following
arguments can be made rigorous by applying, e.\,g., finite differences in time.
We differentiate \eqref{sn3} (formally) with respect to time, obtaining the identity
\begin{align*}
&\partial^2_{tt}\rhon-\Delta\pt\rhon+ \bigl(\beta^\eps)'(\rhon)\,\pt\rhon\\
&=\,(-\pi'(\rhon)+{\cal T}_{\tau_n}(\munmu)\,g''(\rhon))\,\pt\rhon\,+\,\pt{\cal T}_{\tau_n}(\munmu)\,g'(\rhon).
\end{align*}
Multiplication by $\pt\rhon$,
integration over $Q_t$, where $t\in (0,T]$, and (formal) integration by parts, using \eqref{sn4}, yields the 
identity
\begin{align}
\label{est2}
&\frac 12\Bigl(\|\pt\rhon(t)\|_H^2\,+\|\pt\rhong(t)\|^2_{\Hg}\Bigr)
\,+\,\txinto|\nabla\pt\rhon|^2\dx\ds\,+\,\tginto|\nabla_\Gamma\pt\rhong|^2\dg\ds\non\\
&+\txinto\bigl(\beta^\eps)'(\rhon)\,|\pt\rhon|^2\dx\ds\,+\tginto\bigl(\beta^\eps_\Gamma)'(\rhong)\,|\pt\rhong|^2\dg\ds
\non\\
&=\,\frac 12\Bigl(\|\pt\rhon(0)\|^2_H\,+\,\|\pt\rhong(0)\|^2_{\Hg}\Bigr)
\,+\,\txinto (-\pi'(\rhon)+{\cal T}_{\tau_n}(\munmu)\,g''(\rhon))\,|\pt\rhon|^2\dx\dt\non\\[1mm]
&\quad + \txinto\pt{\cal T}_{\tau_n}(\munmu)\,g'(\rhon)\,\pt\rhon\dx\ds
+\tginto(\pt u_\Gamma-\pi_\Gamma'(\rhong)\,\pt\rhong)\,\pt\rhong\dg\ds\,.
\end{align}
By the monotonicity of $\beta^\eps$ and $\beta_\Gamma^\eps$, all of the terms on the left-hand side are 
nonnegative, and, invoking {\bf (A2)}, \eqref{T2}, the boundedness 
of $\pi',\pi'_\Gamma,g',g'',\munmu$, and Young's inequality,
we readily conclude from \eqref{est1} \gianni{and the first \eqref{regtau} with $k=n-1$}
that the three integrals on the right-hand side are finite.

Moreover, by {\bf (A1)}, {\bf (A4)}, {\bf (A5)} and {\bf (A6)},
we have $\,\mu_0\,g'(\rho_0)+\Delta\rho_0-\pi(\rho_0)\in H$, and it follows from \eqref{domeps}
that, for every $\eps>0$,
\beq
\|\beta^\eps(\rho_0)\|_H\,\le\,\|\beta^\circ(\rho)\|_H \quad\mbox{and}\quad
\|\beta^\eps_\Gamma(\rho_{0_\Gamma})\|_{\Hg}\,\le\,\|\beta_\Gamma^\circ(\rho_{0_\Gamma})\|_{\Hg}\,.
\eeq 
We thus can infer from {\bf (A6)} that 
\beq
\label{K1}
\|\pt\rhon(0)\|_H\,=\,\|\mu_0\,g'(\rho_0)+\Delta\rho_0-\pi(\rho_0)-\beta^\eps(\rho_0)\|_H\,\le\,K_1.
\eeq
Likewise, using {\bf (A1)}, {\bf (A2)}, {\bf (A4)}, {\bf (A6)}, and \eqref{domeps}, we conclude
that
\beq\label{K2}
\|\pt\rhong(0)\|_{\Hg}\,=\,\|\ug(0)-\pn\rho_0+\Delta_\Gamma\rho_{0_\Gamma}-\pig(\rho_{0_\Gamma})
-\beta_\Gamma^\eps(\rho_{0_\Gamma})\|_{\Hg}\,\le\,K_2.
\eeq
In conclusion, the expression on the \rhs~of \eqref{est2} is finite, which proves the additional
regularity \eqref{extrareg}.

\vspace{2mm}
Now that the existence of a unique solution $\,\rhon\,$ to the system \eqref{sn3}--\eqref{sn5}
with the required regularity is established, we substitute it in \eqref{sn1} and study the
resulting initial-boundary value problem \eqref{sn1}--\eqref{sn2}. Recalling the continuity of
the embedding
\beq\label{embed}
\bigl(L^\infty(0,t;H)\cap L^2(0,t;V)\bigr)\,\subset\,\bigl(L^{10/3}(Q_t)\cap 
L^{7/3}(0,t;L^{14/3}(\oma))\bigr) \quad\forall\,t\in (0,T],
\eeq
we have that
\beq\label{rhont}
\pt\rhon\in L^{10/3}(Q_{t_n})\cap 
L^{7/3}(0,t_n;L^{14/3}(\oma)).
\eeq
We now recall the form \eqref{extendg} of the extension of $g$ to the whole real line. Using
this, we may argue as in \cite[pp.~7953--7956]{CGS3} to conclude that, with $\rhon$ fixed,
the problem \eqref{sn1}--\eqref{sn2} has a unique nonnegative solution $\mun\in H^1(I_n;H)
\cap L^\infty(I_n;V)\cap L^2(I_n;W^{2,3/2}(\oma))$. But since $\pt\rhon\in L^2(I_n;V)$, it
follows from the continuity of the embedding $V\subset L^6(\oma)$ that $g'(\rhon)\,\pt\rhon\,
\mu_n\in L^2(I_n;H)$, and comparison in \eqref{sn1} shows that $\Delta\mun\in L^2(I_n;H)$,
whence, by standard elliptic estimates, also $\mun\in L^2(I_n;W)$. The continuity of the
embedding $(H^1(I_n;H)\cap L^2(I_n;\Hdue))\subset C^0(\bar I_n;V)$ yields that then also
$\mun\in C^0(\bar I_n;V)$. 

Finally, we have $\pt\rhon\in L^{7/3}(I_n;L^{14/3}(\oma))$, by
\eqref{rhont}. Therefore, we may repeat the argument from the proof 
of \cite[Thm.~2.3]{CGPS3}, which was based on this regularity,
to conclude that $\mun\in L^\infty(\oma\times I_n)$. Here, we remark that the
quoted proof was performed only for the special case $g(r)=r$; however, the argument
extends with only minor modifications to the present situation (see also the
analogous \cite[Thm.~3.7]{CGPS5}). With this, the above claim is 
completely proved.

\vspace{5mm}\noindent
\underline{\sc Step 2:} \quad Let $\eps>0$ and $N\in\nz$ be arbitrary but fixed, and let
$\bigl(\mu^\eps_{\tau_N},\rho^\eps_{\tau_N},\rho^\eps_{\tau_{N_\Gamma}}\bigr)$ be
the unique solution to $(P^\eps_{\tau_N})$ having the regularity properties 
\eqref{regepstau}, which was constructed in Step~1. We are now going to
show a number of a priori estimates for these solutions which are uniform with respect to
$N$ and $\eps>0$. In performing these estimates, for the sake of a better readability \pier{we will}
omit both the subscript $\tau_N$ and the superscript $\eps$, writing them only at the
end of each estimate. We also recall that $K_i$, $i\in\nz$, denote positive constants which
are independent of both $\eps$ and $N$. We remark that the following chain of estimates
combines ideas from \cite{CGPS5} and \cite{Calcol}, where in \cite{CGPS5} a more general
version of Eq.~\eqref{ss1} was investigated. Since, however, in \cite{CGPS5} the simpler case of a zero
Neumann boundary condition for $\,\rho\,$ was assumed in place of the dynamic boundary condition
\eqref{ss5}, we have chosen to include these estimates for the reader's convenience.   

\vspace{2mm}\noindent
\underline{\sc First estimate:} 
\par\nobreak
\vspace{1mm}\noindent
First observe that $\,\,\pt\bigl(\bigl(\mbox{$\frac 12$}+g(\rho)\bigr)\mu^2\bigr)
=(1+2g(\rho))\,\mu\,\pt\mu+g'(\rho)\,\pt\rho\,\mu^2$. Hence, multiplying \eqref{etau1} by
$\mu$, and integrating over $Q_t$, where $t\in (0,T]$, we find the identity
\begin{align*}
\xinto\bigl(\mbox{$\frac 12$}+g(\rho(t))\bigr)\,\mu^2(t)\dx\,+\txinto|\nabla\mu|^2\dx\ds
\,=\,\xinto\bigl(\mbox{$\frac 12$}+g(\rho_0)\bigr)\,\mu_0^2(t)\dx,
\end{align*}
whence, using {\bf (A1)}, we easily conclude that
\beq\label{K3}
\left\|\mu^\eps_{\tau_N}\right\|_{L^\infty(0,T;H)\cap L^2(0,T;V)}\,\le\,K_3,
\eeq
which entails that also
\beq\label{K4}
\left\|{\cal T}_{\tau_N}(\mu^\eps_{\tau_N})\right\|_{L^\infty(0,T;H)\cap
L^2(0,T;V)}\,\le\,K_4.
\eeq

\vspace{2mm}\noindent
\underline{\sc Second estimate:} 
\par\nobreak
\vspace{1mm}\noindent
Next, we add $\rho$ on both sides of \eqref{etau3}, multiply by $\pt\rho$, and integrate over $Q_t$, where $t\in (0,T]$.
Adding $\rho_\Gamma$ on both sides of \eqref{etau4}, we then obtain the identity
\begin{align}
&\txinto|\pt\rho|^2\dx\ds\,+\,\frac 1 2\,\|\rho(t)\|_V^2
+\xinto\pier{\widehat\beta{}}^\eps(\rho(t))\dx\non\\[1mm]   
&+\tginto|\pt\rg|^2\dg\ds\,+\,\frac 12\,\|\rg(t)\|^2_{\Vg}\,+\ginto\pier{\widehat\beta{}}_\Gamma^\eps(\rg(t))\dg
\non\\[1mm]
&=\,\frac 1 2\,\|\rho_0\|_V^2\,+\xinto\pier{\widehat\beta{}}^\eps(\rho_0)\dx
\,+\,\frac 12\,\|\rho_{0_\Gamma}\|_{\Vg}^2\,+\ginto\pier{\widehat\beta{}}_\Gamma^\eps(\rho_{0_\Gamma})\dg
\non\\[1mm]
&\,\,+\txinto\pt\rho\,({\cal T}_{\tau_N}(\mu)\,g'(\rho)-\pi(\rho)+\rho)\dx\ds\,+\tginto\pt\rg\,
(\ug-\pig(\rg)+\rg)\dg\ds\,.
\end{align}
From \eqref{antider2}, we obtain 
that all of the expressions on the \lhs~are nonnegative, 
and it follows from {\bf (A1)}, \eqref{antider2}, and {\bf (A6)},
that the first four summands on the \rhs~are bounded by a constant that neither depends on $N\in\nz$ nor
on $\eps>0$. The last two integrals on the \rhs, which we denote by $I_1$ and $I_2$, respectively,
can be estimated as follows: we employ \eqref{extendg}, \eqref{K4}, {\bf (A4)}, and Young's inequality, 
to obtain that
\begin{equation}\label{K5}
I_1\,\le\,\frac 12\txinto|\pt\rho|^2\dx\ds\,+\,K_5\Bigl(1\,+\txinto|\rho|^2\dx\ds\Bigr),
\end{equation}
and invoke {\bf (A2)}, {\bf (A4)}, and Young's inequality, to see that
\beq\label{K6}
I_2\,\le\,\frac 12\tginto|\pt\rg|^2\dg\ds\,+\,K_6\Bigl(1\,+\tginto|\rg|^2\dg\ds\Bigr).
\eeq
We thus can infer from Gronwall's lemma that
\beq\label{K7}
\left\|\rho^\eps_{\tau_N}\right\|_{H^1(0,T;H)\cap L^\infty(0,T;V)}\,+\,\left\|\rho^\eps_{\tau_{N_\Gamma}}
\right\|_{H^1(0,T;\Hg)\cap L^\infty(0,T;\Vg)}\,\le\,K_7.
\eeq         

\vspace{2mm}\noindent
\underline{\sc Third estimate:}
\par\nobreak
\vspace{1mm}\noindent
Next, we multiply \eqref{etau3} by $\,\beta^\eps(\rho)\,$ and integrate over $\,Q_t$, where $0<t\le T$.
\gianni{Then, we integrate by parts and use \eqref{etau4}--\eqref{etau5}}
to obtain the identity
\begin{align}
\label{est31}
&\xinto\pier{\widehat\beta{}}^\eps(\rho(t))\dx\,+\ginto\pier{\widehat\beta{}}^\eps(\rg(t))\dg \,+\txinto|\beta^\eps(\rho)|^2\dx\ds\non\\
&+\txinto(\beta^\eps)'(\rho)\,|\nabla\rho|^2\dx\ds\,+\tginto(\beta^\eps)'(\rg)\,|\ng\rg|^2\dg\ds\non\\
&=\,\xinto\pier{\widehat\beta{}}^\eps(\rho_0)\dx\,+\ginto\pier{\widehat\beta{}}^\eps(\rho_{0_\Gamma})\dg\,+
\txinto({\cal T}_{\tau_N}(\mu)\,g'(\rho)-\pi(\rho))\,\beta^\eps(\rho)\dx\ds\non\\
&\quad+\tginto(\ug-\pi_\Gamma(\rg))\,\beta^\eps(\rg)\dg\ds\,-\tginto\beta^\eps_\Gamma(\rg)
\,\beta^\eps(\rg)\dg\ds\,.
\end{align}
In view of \eqref{antider2}, and as $\beta^\eps$ is increasing, all of the terms on the \lhs~are nonnegative. \pier{About the last term, we exploit \eqref{domination3} to deduce that
\beq
\label{pier1}
-\tginto\beta^\eps_\Gamma(\rg) \,\beta^\eps(\rg)\dg\ds \leq
-\frac{1}{2\eta} \tginto  |\beta^\eps(\rg)|^2 \dg\ds
+ \frac{T\,C_\Gamma^2}{2\eta} \ginto 1 \dg.
\eeq 
}%
Moreover, by \eqref{antider2} and {\bf (A6)},
\beq\label{K8}
\xinto\pier{\widehat\beta{}}^\eps(\rho_0)\dx\,+ \ginto\pier{\widehat\beta{}}^\eps(\rho_{0_\Gamma})\dg\,\le\xinto\pier{\widehat\beta{}}(\rho_0)\dx
\,+\ginto\pier{\widehat\beta{}}(\rho_{0_\Gamma})\dg\,\le\,K_8.
\eeq
\gianni{Furthermore}, thanks to the previous estimates and Young's inequality,
\beq\label{K9}
\txinto({\cal T}_{\tau_N}(\mu)\,g'(\rho)-\pi(\rho))\,\beta^\eps(\rho)\dx\ds\,
\le\,\frac 12\txinto|\beta^\eps(\rho)|^2\dx\ds
\,+\,K_9,
\eeq
as well as 
\beq\label{K10}
\tginto(\ug-\pi_\Gamma(\rg))\,\beta^\eps(\rg)\dg\ds\,\le\,\frac 1{4\eta}\tginto|\beta^\eps(\rg)|^2\dg\ds\,+\,K_{10},
\eeq
with the constant $\eta>0$ introduced in {\bf (A7)}. Hence, 
\pier{we can infer from \eqref{est31}--\eqref{K10}} that
\beq\label{K11}
\left\|\beta^\eps(\rho^\eps_{\tau_N})\right\|_{L^2(0,T;H)}\,\le\,K_{11},
\eeq
whence, by comparison in \eqref{etau3},
\beq
\label{K12}
\left\|\Delta\rho_{\tau_N}^\eps\right\|_{L^2(0,T;H)}\,\le\,K_{12}.
\eeq 

\vspace{2mm}\noindent
\underline{\sc Fourth estimate:}
\par\nobreak
\vspace{1mm}\noindent
We now draw some consequences from \eqref{K12}. First, we invoke \cite[Theorem~3.1, p.~1.79]{brez}, 
which yields the estimate
\begin{equation}
\label{K13}
\int_0^T\|\rho(t)\|^2_{H^{3/2}(\Omega)}\dt\,
\leq K_{13} \int_0^T(\|\Delta \rho(t)\|^2_{H}+\|\rho_\Gamma(t)\|^2_{V_\Gamma})\dt,
\end{equation}
whence we obtain that (cf.~\eqref{K7}) 
\begin{equation}
\label{K14}
\|\rho\|_{L^2(0,T;H^{3/2}(\Omega))}\,\le\,K_{14}.
\end{equation}
By virtue of \eqref{K12} and of
\cite[Theorem~2.27, p.~1.64]{brez}, we can also conclude that
\begin{equation}\label{K15}
\left\|\pn \rho\right\|_{L^2(0,T; H_\Gamma) } \,\le\,K_{15}.
\end{equation}

We now aim to find an estimate for $\,\beta^\eps_\Gamma(\rg)$. To this end,
 we multiply \eqref{etau4} by $\beta_\Gamma^\varepsilon(\rg)$ and integrate over $\Gamma\times(0,t)$,
where $0<t\le T$. We obtain the identity
\begin{align}\label{esti41}
&\ginto\pier{\widehat\beta{}}^\eps_\Gamma(\rg(t))\dg  \,+ \tginto(\beta^\eps_\Gamma)'(\rg)\,|\ng\rg|^2\dg\ds\,+\,\tginto|\beta_\Gamma^\eps(\rg)|^2\dg\ds\non\\
&=\ginto\pier{\widehat\beta{}}^\eps_\Gamma(\rho_{0_\Gamma})\dg\,+\tginto(\ug-\pn\rho-\pi_\Gamma(\rg))\,
\beta^\eps_\Gamma(\rg)\dg\ds\,.
\end{align}
By the monotonicity of $\,\beta^\eps_\Gamma\,$ and  \eqref{antider2}, the first two integrals
 on the left-hand side of \eqref{esti41} are nonnegative, while, thanks to
{\bf (A6)},
\begin{equation}\label{K16}
\int_\Gamma \pier{\widehat\beta{}}_\Gamma^\varepsilon(\rho_{0_\Gamma})\dg \,\le\ginto
\pier{\widehat\beta{}}_{\Gamma}(\rho_{0_\Gamma})\dg\,\le\,K_{16}.
\end{equation}
For the last integral on the right-hand side, we invoke {\bf (A2)}, {\bf (A4)}, \eqref{K7},
\eqref{K15}, and Young's inequality,
 to infer that
\begin{equation}\label{K17}
\tginto ( u_\Gamma-\pn\rho-\pi_\Gamma(\rho_\Gamma))\,\beta_\Gamma^\varepsilon(\rg)\dg\ds\,\le
\frac{1}{2}\,\tginto|\beta_\Gamma^\varepsilon(\rg)|^2\dg\ds\,+\,K_{17},
\end{equation}
whence we conclude that
\beq\label{K18}
\left\|\beta^\eps_\Gamma(\rg)\right\|_{L^2(0,T;\Hg)}\,\le\,K_{18}.
\eeq
Comparison in \eqref{etau4}, using the previously shown estimates, then \pier{yields}
\beq\label{K19}
\|\Dg\rg\|_{L^2(0,T;\Hg)}\,\le\,K_{19},
\eeq
and it follows from the boundary version of the elliptic estimates that
\beq\label{K20}
\|\rg\|_{L^2(0,T;H^2(\Gamma))}\,\le\,K_{20},
\eeq
whence, employing \eqref{K12} and standard elliptic estimates, we conclude that also
\beq\label{K21}
\|\rho\|_{L^2(0,T;\Hdue)}\,\le\,K_{21}.
\eeq
In conclusion, we have \pier{inferred} the estimate
\beq\label{K22} 
\left\|\beta^\eps_\Gamma(\rho^\eps_{\tau_{N_\Gamma}})\right\|_{L^2(0,T;\Hg)}\,+\,
\left\|\rho^\eps_{\tau_N}\right\|_{L^2(0,T;\Hdue)}\,+\,\left\|\rho^\eps_{\tau_{N_\Gamma}}\right\|
_{L^2(0,T;H^2(\Gamma))}\,\le\,K_{22}.
\eeq

\vspace{2mm}\noindent
\underline{\sc Fifth estimate:}
\par\nobreak
\vspace{1mm}\noindent
In this step of the proof, we partly repeat a formal argument that was already used in Step~1, noting again
that it can be made rigorous by using finite differences in time. Namely, we differentiate \eqref{etau3}
formally with respect to time, multiply the resulting identity by $\,\pt\rho$, and integrate over $Q_t$,
where $0<t\le T$, and (formally) by parts. As in Step~1 (see the argumentation between Eqs.~\eqref{est2} and
\eqref{K2}), we then arrive at an inequality of the form
\begin{align}\label{K23}
&\frac 12\left(\|\pt\rho(t)\|_H^2+\|\pt\rg(t)\|_{\Hg}^2\right)\,+\txinto|\nabla\pt\rho|^2\dx\ds
\,+\tginto|\ng\pt\rg|^2\dg\ds\non\\
&\le\,K_{23}\,+\,\sum_{j=1}^4 I_j,
\end{align}    
where the expressions $\,I_j$, $1\le j\le 4$, will be specified and estimated below. At first, recall that
$\,\pi',\pi'_\Gamma,g',g''$\, are bounded. Hence, owing to \eqref{K7},
\pier{we have that}
\beq\label{K24}
I_1:=-\txinto \pi'(\rho)\,|\pt\rho|^2\dx\ds\,\le\,K_{24},
\eeq  
\pier{and}, by also using {\bf (A2)} and Young's inequality, \pier{we infer that}
\beq\label{K25}
I_4:=\tginto(\pt\ug-\pig'(\rg)\,\pt\rg)\,\pt\rg\dg\ds\,\le\,K_{25}.
\eeq
In addition, invoking  
\eqref{K4} and the compactness inequality \eqref{compact}, as well as H\"older's 
inequality, we find that
\begin{align}
\label{K26}
&I_2:=\txinto{\cal T}_{\tau_N}(\mu)\,g''(\rho)\,|\pt\rho|^2\dx\ds\,\le\,
K_{26}\int_0^t\left\|{\cal T}_{\tau_N}(\mu(s))\right\|_2\,\|\pt\rho(s)\|_4^2\ds\non\\
&\hspace*{7mm}\le K_{27}\int_0^t\|\pt(\rho(s)\|_4^2\ds\,\le\,\frac 16\int_0^t\|\nabla\pt\rho(s)\|_H^2\ds\,+\,K_{28}.
\end{align} 

The estimation of the remaining term
$$
I_3:=\txinto\pt{\cal T}_{\tau_N}(\mu)\,g'(\rho)\,\pt\rho\dx\ds
$$
is more delicate. To this end, we note that \eqref{etau1} implies the identity
\beq\label{muet}
\pt\mu\,=\,\frac 1{1+2g(\rho)}\bigl(\Delta\mu-\mu\,g'(\rho)\,\pt\rho\bigr),
\eeq
where, thanks to \eqref{extendg},
$\,\,1/(1+2g(\rho))\le 3$. 
First, recall that ${\cal T}_{\tau_N}(\mu)$ is constant with respect to time on the interval~$(0,\tau)$.
We therefore have, with obvious notation:
\begin{align}
& \txinto g'(\rho)\,\pt{\cal T}_{\tau_N}(\mu) \, \pt\rho\dx\ds
  \,= \int_0^{t-\tau_N}\!\!\!\xinto \pt\mu(s) \, g'(\rho(s+\tau_N))\, \pt\rho(s+\tau_N)\dx \ds
  \non
  \\
& =\, \int_0^{t-\tau_N}\!\!\!\xinto \frac 1{\coeffs} \,
  \bigl[ \Delta \mu(s) - \mu(s) \,g'(\rho(s))\, \pt\rho(s) \bigr]\,
  \pt g(\rho(s+\tau_N)) \dx\ds
  \non
  \\[1mm]
&=\,-\int_0^{t-\tau_N}\!\!\!\xinto \nabla\mu(s) \cdot \nabla \Bigl(
  \frac {\pt g(\rho(s+\tau_N))} {\coeffs} \Bigr)\dx\ds
  \non
  \\[1mm]
&  \quad\,  - \int_0^{t-\tau_N}\!\!\!\xinto \frac {g'(\rho(s)) \, g'(\rho(s+\tau_N)) }{\coeffs} \, \mu(s)\,\pt\rho(s) \,\pt\rho(s+\tau_N)\dx \ds\,. 
  \label{est52}
\end{align}
We treat the last two integrals separately, invoking our structural assumptions and the estimates established above.
We have 
\begin{eqnarray}\label{K29}
  && -\int_0^{t-\tau_N}\!\!\!\xinto \nabla\mu(s) \cdot \nabla \Bigl(\frac 
  {\pt g(\rho(s+\tau_N))}  {\coeffs}
  \Bigr) \dx\ds \non
  \\
  &&= -\int_0^{t-\tau_N}\!\!\!\xinto \nabla\mu(s) \cdot\nabla \Bigl( 
  \frac{g'(\rho(s+\tau_N)) \,\pt \rho(s+\tau_N)}{\coeffs}\Bigr)\dx \ds
  \non\\
  && \leq \,K_{29}\int_0^{t-\tau_N}\!\!\!\xinto |\nabla\mu(s)| \, |\nabla\pt\rho(s+\tau_N)|\dx \ds
  \non 
   \\
  && \quad +\,K_{30}\int_0^{t-\tau_N}\!\!\!\xinto |\nabla\mu(s)| \, |\nabla\rho(s)| 
  \, |\pt\rho(s+\tau_N)|\dx \ds
  \non  \\
  && \quad +\,K_{31}\int_0^{t-\tau_N}\!\!\!\xinto |\nabla\mu(s)| \, 
  |\nabla\rho(s +\tau_N)| \, |\pt\rho(s+\tau_N)|
  \dx \ds, 
\end{eqnarray}
where, owing to \eqref{K3} and Young's inequality,
\beq\label{K32}
\pier{K_{29}}\int_0^{t-\tau_N}\!\!\!\xinto |\nabla\mu(s)| \, |\nabla\pt\rho(s+\tau_N)|\dx\ds
\,\le\,
\frac 16 \txinto  |\nabla\pt\rho|^2 \dx\ds\,+\,K_{32}.
\eeq  

On the other hand, we also have that
\beqa\label{K33}
  && \pier{K_{30}}\int_0^{t-\tau_N}\!\!\!\xinto |\nabla\mu(s)| \, |\nabla\rho(s)| \, |\pt\rho(s+\tau_N)| \dx\ds
  \non  \\
  && \leq \pier{K_{30}}\int_0^{t-\tau_N} \|\nabla\mu(s)\|_2\,\|\nabla\rho(s)\|_4 \,\|\pt\rho(s+\tau_N)\|_4 \ds
  \non  \\
  && \leq \,\frac 16\int_0^t\|\pt\rho(s)\|_V^2 \ds \,+\,K_{33} \int_0^{t-\tau_N} \|\nabla\mu(s)\|_H^2\,
   \|\nabla\rho(s)\|_V^2 \ds
  \non \\
  && \leq \,\frac 16 \txinto |\nabla\pt\rho|^2\dx\ds\,
  + K_{34}\txinto |\pt\rho|^2\dx\ds \,+
   K_{35} \int_0^{t-\tau_N}  \|\nabla\mu(s)\|_H^2 \|\nabla\rho(s)\|_V^2 \ds
  \non
  \\
  && \leq \frac 16 \txinto |\nabla\pt\rho|^2\dx\ds \,+\,K_{36}\,  
  + \,K_{37} \int_0^{t-\tau_N} \|\nabla\mu(s)\|_H^2 \,\|\nabla\rho(s)\|_V^2 \ds.
\eeqa

The last integral cannot be controlled in this form. We thus try to estimate
the expression $\,\|\nabla\rho(s)\|_V^2	\,$ in terms of the expressions
\,$\,\|\pt\rho(s)\|_H^2\,$ and $\,\|\pt\rg(s)\|^2_{\Hg}\,$ which can be handled 
using the first summand on the \lhs~of \eqref{K23}. To this end, we use the 
regularity theory for linear elliptic equations and \eqref{K7} to deduce that 
\beq\label{K38}
\|\nabla\rho(s)\|_V^2\,\le\, K_{38}(\|\rho(s)\|_V^2\,+\,\|\Delta\rho(s)\|_H^2)
\,\le\,K_{39}\left(1\,+\,\|\Delta\rho(s)\|_H^2\right).
\eeq
We now multiply, just as in the third estimate above, \eqref{etau3} by $\,\beta^\eps(\rho(s))$,
but this time we only integrate over $\oma$ to obtain the identity (compare \pier{with} \eqref{est31})
\begin{align}
\label{est53}
&\|\beta^\eps(\rho(s))\|_H^2\,+\xinto\bigl(\beta^\eps\bigr)'(\rho(s))\,|\nabla\rho(s)|^2\dx
 \,+\,\ginto \bigl(\beta^\eps\bigr)'(\rg(s))\,|\ng\rg(s)|^2\dg\non\\
  &=\,\xinto\beta^\eps(\rho(s))\left(-\pt\rho(s)-\pi(\rho(s))+{\cal T}_{\tau_N}(\mu(s))\,
g'(\rho(s))\right)\dx\non\\
&\quad +\ginto\beta^\eps(\rg(s))\left(\ug(s)-\pig(\rg(s)) -\pt\rg(s)\right)\dg
 \pier{{}- \ginto \beta^\eps_\Gamma(\rg(s))\beta^\eps(\rg(s)) \,\dg}
\,.  
\end{align}
The terms on the \lhs~are nonnegative, and analogous reasoning as in the third estimate,
using the already proven \gianni{bounds}, the general assumptions, \eqref{K4}, \pier{\eqref{domination3}} and Young's inequality, shows that
\beq\label{K40}
\|\beta^\eps(\rho(s))\|^2_H\,\le\,K_{40}\left(1+\|\pt\rho(s)\|_H^2+\|\pt\rg(s)\|^2_{\Hg}\right)
\quad\mbox{for a.\,e. $s\in (0,t)$,}
\eeq
whence, by comparison in \eqref{etau3},
\beq\label{K41}
\|\Delta\rho(s)\|_H^2\,\le\,K_{41}\left(1+\|\pt\rho(s)\|_H^2+\|\pt\rg(s)\|^2_{\Hg}\right)
\quad\mbox{for a.\,e. $s\in (0,t)$.}
\eeq 
Combining the estimates \eqref{K33}--\eqref{K41}, we have thus shown that
\begin{align}
&\int_0^{t-\tau_N}\!\!\!\xinto|\nabla\mu(s)| \, |\nabla\rho(s)| \, |\pt\rho(s+\tau_N)|\dx \ds
  \non  \\
&\leq\,\frac 16 \txinto |\nabla\pt\rho|^2 \dx\ds\,+\,
K_{42}\Bigl(1+\int_0^t\|\nabla\mu(s)\|_H^2 \left(\|\pt\rho(s)\|_H^2 
\,+\,\|\pt\rg(s)\|^2_{\Hg}\right)\ds\Bigr).
\label{K42}
\end{align}
By similar reasoning, we also obtain that
\begin{align}
\pier{K_{31}}\int_0^{t-\tau_N}\!\!\!\xinto|\nabla\mu(s)| \, |\nabla\rho(s+\tau_N)| \, |\pt\rho(s+\tau_N)|
\dx \ds  \,\le\,\frac 16 \txinto |\nabla\pt\rho|^2\dx\ds\non\\ 
+ \,K_{43}\Bigl(1+\int_0^t\|\nabla{\cal T}_{\tau_N}(\mu(s))\|_H^2 \left(\|\pt\rho(s)\|_H^2 
\,+\,\|\pt\rg(s)\|^2_{\Hg}\right)\ds\Bigr), 
  \label{K43}
\end{align}
where \,$\nabla{\cal T}_{\tau_N}(\mu(s))=\nabla\mu_0\,$ for $0\le s\le \tau_N$.
Hence, combining the \pier{inequalities \eqref{K29}--\eqref{K32} and \eqref{K42}--\eqref{K43}}, we have proved the estimate 
\begin{align}
&-\int_0^{t-\tau_N}\!\!\!\xinto \nabla\mu(s) \cdot \nabla \Bigl(\frac 
  {\pt g(\rho(s+\tau_N))} {\coeffs}\Bigr)\dx \ds\, \leq \,\frac 12\txinto |\nabla\pt\rho|^2\dx\ds
  \non  \\[1mm]
&+\,K_{44}\Bigl(1\,+\int_0^t \left(\|\nabla\mu(s)\|_H^2 \,+\,\|\nabla{\cal T}_{\tau_N}(\mu(s))\|_H^2
   \right)\left(\|\pt\rho(s)\|_H^2\,+\,\|\pt\rg(s)\|_{\Hg}^2\right) \ds\Bigr)\,.
  \label{K44}
\end{align}

 Let us now come to the second term on the \rhs~of \eqref{est52}. We have, 
 owing to H\"older's and Young's inequalities,
\begin{align}
&-\int_0^{t-\tau_N}\!\!\!\xinto \frac {g'(\rho(s))\,g'(\rho(s+\tau_N))}{\coeffs} \, \mu(s)
\, \pt\rho(s)\, \pt\rho(s+\tau_N)\dx \ds 
  \non  \\
& \leq\,K_{45}\int_0^{\gianni{t-\tau_N}}\|\mu(s)\|_3\,\|\pt\rho(s+\tau_N)\|_6\, \|\pt\rho(s)\|_2 \ds
  \non  \\
& \leq \frac 16\int_0^{t-\tau_N} \|\pt\rho(s+\tau_N)\|_V^2 \ds\,+ 
  \, K_{46} \int_0^t \|\mu(s)\|_3^2 \,\|\pt\rho(s)\|_H^2\ds
  \non  \\
& \leq \frac 16 \int_0^t\|\nabla\pt\rho(s)\|_H^2 \ds\, + \,K_{47}\Bigl(1+
  \int_0^t \|\mu(s)\|_3^2\,\|\pt\rho(s)\|_H^2 \ds\Bigr) .
  \label{K45}
\end{align}
Therefore, due to \eqref{K44} and \eqref{K45}, we can infer from \gianni{\eqref{est52}} that
\begin{align}
 &\txinto g'(\rho)\,\pt{\cal T}_{\tau_N}(\mu) \, \pt\rho\dx\ds
  \,\leq \,\frac 23\txinto |\nabla\pt\rho|^2\dx\ds\,+\,K_{48}\non\\
 &+\,K_{49}\int_0^t
   \gianni{\Phi(s)}
   \left(\|\pt\rho(s)\|_H^2
    +\|\pt\rg(s)\|_{\Hg}^2
   \right) \ds,
    \label{K48}
\end{align}
\gianni{where $\Phi:(0,T)\to\rz$ is given by 
$\Phi(s):=\|\mu(s)\|_3^2+\|\nabla\mu(s)\|_H^2+\|\nabla{\cal T}_{\tau_N}(\mu(s))\|_H^2\,$
for $s\in(0,T)$.
Since it is known that $\Phi$ is} bounded in $L^1(0,T)$ uniformly with respect to
$N\in\nz$ and $\eps>0$, summarizing the estimates \eqref{K23}--\eqref{K26}, \eqref{K48}, 
and invoking Gronwall's lemma, we have thus finally shown that
\beq
 \left \|\pt\rho^\eps_{\tau_N}\right\|_{L^\infty(0,T; H)\cap L^2(0,T;V)}
 \,+\,\left\|\pt\rho^\eps_{\tau_{N_\Gamma}}\right\|_{L^\infty(0,T;\Hg)
 \cap L^2(0,T;\Vg)}\, \leq\,K_{50}.
  \label{K50}
\eeq

Let us draw some consequences from this estimate. First note that \eqref{K40},
\eqref{K41} and \eqref{K50} entail that
\beq\label{K51}
\left\|\beta^\eps(\rho^\eps_{\tau_N})\right\|_{L^\infty(0,T;H)}\,+\,
\left\|\Delta\rho^\eps_{\tau_N}\right\|_{L^\infty(0,T;H)}\,\le\,K_{51}. 
\eeq
Next, we recall \eqref{K7} and the fact that $\,\ug\in C^0([0,T];\Hg)$, by 
assumption {\bf (A2)}. We therefore can, for almost every fixed $t\in (0,T)$,
follow exactly the same chain of arguments as in the fourth estimate, but
this time without integrating over time. In doing this, we establish consecutively
bounds resembling the estimates \eqref{K13}--\eqref{K22},
eventually arriving at the estimate
\begin{align}
\label{K52}
\left\|\beta^\eps_\Gamma(\rho^\eps_{\tau_{N_\Gamma}})\right\|_{L^\infty(0,T;\Hg)}\,
+\left\|\rho^\eps_{\tau_N}\right\|_{L^\infty(0,T;\Hdue)}\,+\left\|
\rho^\eps_{\tau_{N_\Gamma}}\right\|_{L^\infty(0,T;H^2(\Gamma))}\,\le\,K_{52}.
\end{align}  
Notice that \eqref{K52} implies, in particular, that
\beq\label{K53}
\left\|\rho^\eps_{\tau_N}\right\|_{L^\infty(Q)}\,+\left\|\rho^\eps_{\tau_{N_\Gamma}}
\right\|_{L^\infty(\Sigma)}\,\le\,K_{53}.
\eeq

\vspace{2mm}\noindent
\underline{\sc Sixth estimate:}
\par\nobreak
\vspace{1mm}\noindent
We now multiply \eqref{etau1} by $\,\pt\mu\,$ and integrate over $Q_t$,
where $0<t\le T$, and by parts. Recalling \eqref{extendg} and {\bf (A1)}, and
invoking H\"older's and Young's inequalities, we obtain that
\begin{align}\label{K54}
&\frac 13\txinto|\pt\mu|^2\dx\ds\,+\,\frac 12\,\|\nabla\mu(t)\|_H^2\,\le\,
\frac 12\,\|\nabla\mu_0\|_H^2\,+\,K_{54}\txinto|\pt\mu|\,|\mu|\,|\pt\rho|\dx\ds
\non\\
&\le\,K_{55}\,+\,K_{56}\int_0^t\|\pt\mu(s)\|_2\,\gianni\|\pt\rho(s)\|_4\,\|\mu(s)\|_4\ds\non\\
&\le\,\gianni{K_{55}\,+\,}\frac 16\txinto|\pt\mu|^2\dx\ds\,+\,K_{57}\int_0^t\|\pt\rho(s)\|_V^2
\,\|\mu(s)\|_V^2\ds
\end{align}
\gianni{and we note that the function $s\mapsto\|\pier{\partial_t}\rho(s)\|_V^2$ 
is bounded in $L^1(0,T)$ by \eqref{K50}. Thus,}
taking \eqref{K3} into account, we \juerg{can infer} from Gronwall's lemma that
\beq\label{K58}
\left\|\mu^\eps_{\tau_N}\right\|_{H^1(0,T;H)\cap L^\infty(0,T;V)}\,\le\,K_{58}.
\eeq
But then, owing to H\"older's inequality and and \eqref{K50},
\begin{align*}
&\|g'(\rho)\,\pt\rho\,\mu\|_{L^2(Q)}^2\,\le\,K_{59}\int_0^T \|\pt\rho(s)\|_4^2\,
\|\mu(s)\|_4^2\ds\non\\
&\le\,K_{60}\,\|\mu\|^2_{L^\infty(0,T;V)}
\,\|\pt\rho\|^2_{L^2(0,T;V)}\,\le\,K_{61},
\end{align*}
whence, by comparison in \eqref{etau1},
$$\|\Delta\mu\|_{L^2(0,T;H)}\,\le\,K_{62}.$$
In view of the boundary condition \eqref{etau2}, we therefore can infer from 
standard elliptic estimates that
\beq\label{K63}
\left\|\mu^\eps_{\tau_N}\right\|_{L^2(0,T;\Hdue)}\,\le\,K_{63}.
\eeq
Finally, we recall \eqref{K50} and \eqref{embed}, which yield that
$\,\,\pt\rho^\eps_{\tau_N}\,\,$ is bounded in the space 
$\,L^{7/3}(0,T;L^{14/3}(\oma))$,
uniformly in $N\in\nz$ and $\eps>0$. Hence, we may again argue as in the 
proof of \cite[Thm.~2.3]{CGPS3} to conclude that
\beq\label{K64}
\left\|\mu^\eps_{\tau_N}\right\|_{L^\infty(Q)}\,\le\,K_{64}.
\eeq

\vspace{2mm}\noindent
\pier{{}\underline{\sc Step 3:} \quad 
In this step} of the proof, we perform the limit process as $N\to\infty$. Owing to
the a priori estimates established above in Step~2, we can select a subsequence,
which is again indexed by $N$, such that, as $N\to\infty$,
\beqa
\label{conN1}
&\rho^\eps_{\tau_N}\,\to\,\rho^\eps&\mbox{weakly-star in $W^{1,\infty}(0,T;H)
\cap H^1(0,T;V)\cap L^\infty(0,T;\Hdue)$,}\quad\quad\\
\label{conN2}
&\rho^\eps_{\tau_{N_\Gamma}}\to\rho_\Gamma^\eps&
\mbox{weakly-star in $W^{1,\infty}(0,T;H_\Gamma)
\cap H^1(0,T;V_\Gamma)\cap L^\infty(0,T;H^2(\Gamma))$,}\quad\quad\\
\label{conN3}
&\mu^\eps_{\tau_N}\to\mu^\eps&\mbox{weakly in $H^1(0,T;H)\cap L^2(0,T;W)$\,
and weakly-star in $\,L^\infty(Q)$},
\eeqa
for suitable limits $\,\mu^\eps,\rho^\eps,\rho^\eps_\Gamma$. By \eqref{conN3},
we have \pier{that} $\,\pn\mu^\eps=0\,$ almost everwhere on \,$\Sigma$, and the continuity
of the embedding $\,H^1(0,T;H)\cap L^2(0,T;\Hdue)\subset C^0([0,T];V)$ implies
that both $\rho^\eps(0)=\rho_0\,$ and $\,\mu^\eps(0)=\mu_0$. Moreover, by virtue
of standard compactness results (see, e.\,g., \cite[Sect.~8,~Cor.~4]{Simon}), we
can without loss of generality assume that\pier{%
\begin{align}\label{conC0}
&\mu^\eps_{\tau_N}\to \mu^\eps\quad\mbox{strongly in $C^0([0,T];H)$},\\
\label{conC1}
&\rho^\eps_{\tau_N}\to\rho^\eps\quad\mbox{strongly in $C^0([0,T];H^s(\oma))\, $ for
$\, 0<s<2$},
\end{align}
}whence we obtain that $\mu^\eps\ge 0$ almost everywhere in $Q$ and 
\beq\label{conC2}
\rho^\eps_{\tau_N}\to\rho^\eps\quad\mbox{strongly in $C^0(\overline{Q})$.}
\eeq
But then $\,\,\rho^\eps_{\tau_N|\pier{\Sigma}}\to\rho^\eps_{|\pier{\Sigma}}\,$ strongly in
$\,C^0(\overline{\Sigma})$, and thus, invoking \eqref{conN2}, we obtain that
$\,\rho^\eps_{|\pier{\Sigma}}=\rho^\eps_\Gamma$. In addition, the Lipschitz continuity 
of the corresponding functions on $\rz$ yields that
\begin{align}
\label{conN4}
&\Phi(\rho^\eps_{\tau_N})\to \Phi(\rho^\eps)\quad\mbox{strongly in 
$\,C^0(\overline{Q})\,$ for $\,\Phi\in\{\beta^\eps, \pi,g,g'\}$},\non\\
&\Phi_\Gamma(\rho^\eps_{\tau_{N_\Gamma}})\to \Phi_\Gamma(\rho^\eps_\Gamma)
\quad\mbox{strongly in 
$\,C^0(\overline{\Sigma})\,$ for $\,\Phi_\Gamma\in\{\beta^\eps_\Gamma, \pi_\Gamma\}$.}
\end{align}
It is therefore easily verified that
\begin{align}\label{conN5}
&\left(1+2g(\rho^\eps_{\tau_N})\right)\pt\mu^\eps_{\tau_N}\to (1+2g(\rho^\eps))
\,\pt\mu^\eps \quad\mbox{weakly in $\,L^2(Q)$},\non\\
&\,\,\mu^\eps_{\tau_N}\,g'(\rho^\eps_{\tau_N})\,\pt\rho^\eps_{\tau_N}
\to \mu^\eps\,g'(\rho^\eps)\,\pt\rho^\eps\quad\mbox{weakly in 
$\,\pier{L^1(Q)}$.}
\end{align}
Finally, in view \pier{of \eqref{K58}}, it is easy to check that 
$$
\|{\cal T}_{\tau_N}(\mu^\eps_{\tau_N})-\mu^\eps_{\tau_N}\|_{\pier{L^2(Q)}}\,\to\,0\quad\mbox{as $\,N\to\infty$},
$$
which, in combination with \pier{\eqref{conC0} and} \eqref{conN4}, implies that
\beq\label{conN6}
{\cal T}_{\tau_N}(\mu^\eps_{\tau_N})\,g'(\rho^\eps_{\tau_N})\,\to\,\mu^\eps\,g'(\rho^\eps)
\quad\mbox{strongly in \,$\pier{L^2(Q)}$}.
\eeq
Hence, letting $N\to\infty$ in the system \eqref{etau1}--\eqref{etau5}, 
we find that the triple $\,(\mu^\eps,\rho^\eps,\rho^\eps_\Gamma)\,$ is a 
solution
to the initial-boundary value problem
\begin{align}
  & \bigl(1 + 2 g(\rho^\eps)\bigr)\, \pt\mu^\eps + g'(\rho^\eps)\,
	\pt\rho^\eps\,\mu^\eps- \Delta\mu^\eps = 0
   \quad \mbox{a.\,e. in $Q$},
  \label{eps1}
  \\
   \label{eps2}
  &\pn\mu^\eps=0\quad\mbox{a.\,e. on $\Sigma$},\\
  \label{eps3}
  & \pt\rho^\eps - \Delta\rho^\eps + \beta^\eps(\rho^\eps)+
	\pi(\rho^\eps) = \mu^\eps\,g'(\rho^\eps)
  \quad \mbox{a.\,e. in $Q$},\\
  \label{eps4}
  &\pn\rho^\eps+\pt\rho^\eps_\Gamma+\beta^\eps_\Gamma(\rho^\eps_\Gamma)+
	\pig(\rho^\eps_\Gamma) -\Dg	\rho^\eps_\Gamma= \ug,
	\quad \rho^\eps_\Gamma=\rho^\eps_{|\Sigma},\quad \mbox{a.\,e. on $\Sigma$},\\
  \label{eps5}
   & \mu^\eps(0) = \mu_0
  \quad\mbox{and}\quad
  \rho^\eps(0) = \rho_0
  \quad \mbox{a.\,e. in $\oma$}, \quad \rho^\eps_\Gamma(0)={\rho_{0}}_{|\Gamma}
	\quad\mbox{a.\,e. on $\Gamma$}.
  \end{align}
Moreover, $\mu^\eps$ is nonnegative almost everywhere in $\,Q$, and,
by virtue of the weak \pier{(and weak star)} sequential lower \pier{semicontinuity} of norms, we
have the estimate
\begin{align}\label{K65}
&\|\rho^\eps\|_{W^{1,\infty}(0,T;H)\cap H^1(0,T;V)\cap L^\infty(0,T;\Hdue)}
\non\\
&+\,\|\rho^\eps_\Gamma\|_{W^{1,\infty}(0,T;\Hg)\cap H^1(0,T;\Vg)
\cap L^\infty(0,T;H^2(\Gamma))}\non\\
&+\,\|\mu^\eps\|_{H^1(0,T;H)\cap L^2(0,T;W)\cap L^\infty(Q)}\non\\
&+\,\|\beta^\eps(\rho^\eps)\|
_{L^\infty(0,T;H)}\,+\,\|\beta^\eps_\Gamma(\rho^\eps_\Gamma)\|
_{L^\infty(0,T;\Hg)}\,\le\,K_{65}.\quad
\end{align}

\vspace{2mm}\noindent
\pier{{}\underline{\sc Step 4:} \quad 
In the final step} of the existence proof, we perform the limit process as $\eps\searrow0$.
Repeating the arguments of Step~3, we can find a sequence $\eps_n\searrow0$ and 
corresponding solution triples $\,(\mu^{\eps_n},\rho^{\eps_n},\rho^{\eps_n}_\Gamma)\,$
to the  system \eqref{eps1}--\eqref{eps5} such that, as $n\to\infty$,
\begin{align}
\label{coneps1}
&\rho^{\eps_n}\,\to\,\rho\quad\mbox{weakly-star in \,\,$W^{1,\infty}(0,T;H)
\cap H^1(0,T;V)\cap L^\infty(0,T;\Hdue)$,}\\
\label{coneps2}
&\rho^{\eps_n}_\Gamma\to\rho_\Gamma\quad
\mbox{weakly-star in \,\,$W^{1,\infty}(0,T;H_\Gamma)
\cap H^1(0,T;V_\Gamma)\cap L^\infty(0,T;H^2(\Gamma))$,}\\
\label{coneps3}
&\mu^{\eps_n}\to\mu\quad\mbox{weakly in $\,H^1(0,T;H)\cap L^2(0,T;W)$\,
and weakly-star in $\,L^\infty(Q)$},\\
\label{coneps4}
&\beta^{\eps_n}(\rho^{\eps_n})\to\xi\quad\mbox{weakly-star in $L^\infty(0,T;H)$,}\\
\label{coneps5}
&\beta_\Gamma^{\eps_n}(\rho^{\eps_n}_\Gamma)\to\xi_\Gamma\quad
\mbox{weakly-star in $L^\infty(0,T;\Hg)$},
\end{align}
for suitable limits $\,\mu,\rho,\rho_\Gamma,\xi,\xi_\Gamma$. As in Step~3, we may 
assume that
\beq\label{conC3}
\mu^{\eps_n}\to \mu\mbox{\, strongly in $C^0([0,T];H)$},
\quad \rho^{\eps_n}\to\rho\mbox{\, strongly in $C^0(\overline{Q})$}, 
\eeq
whence we again conclude that $\,\rho_\Gamma=\rho_{|\Sigma}$, 
$\rho(0)=\rho_0$, and $\mu(0)=\mu_0$. Moreover, \eqref{coneps3} implies that $\pn\mu=0$
almost everywhere on $\Sigma$. We also have that
\begin{align}
&\,\,\Phi(\rho^{\eps_n})\to \Phi(\rho)\quad\mbox{strongly in 
$\,C^0(\overline{Q})\,$ for $\,\Phi\in\{\pi,g,g'\}$},\non\\
&\,\,\pig(\rho^{\eps_n}_\Gamma)\to \pig(\rho_\Gamma)
\quad\mbox{strongly in 
$\,C^0(\overline{\Sigma})\,$,}\non\\[1mm]
&\left(1+2g(\rho^{\eps_n})\right)\pt\mu^{\eps_n}\to (1+2g(\rho))
\,\pt\mu \quad\mbox{weakly in $\,L^2(Q)$},\non\\[1mm]
& \,\,\mu^{\eps_n}\,g'(\rho^{\eps_n})\,\pt\rho^{\eps_n}
\to \mu\,g'(\rho)\,\pt\rho\quad\mbox{weakly in $\,\pier{L^1(Q)}$},\non\\[1mm]
&\,\,\mu^{\eps_n}\,g'(\rho^{\eps_n})\to \mu	\,g'(\rho)\quad\mbox{strongly in $\,\pier{L^2(Q)}.$}\non 
\end{align}
Moreover, by a standard argument in the theory of maximal monotone operators
(see, e.\,g.\pier{,} \cite[Lemma~2.3, p.~38]{Barbu}), we infer that also 
$\,\rho\in D(\beta)\,$ and $\,\xi\in\beta(\rho)\,$ almost
everywhere in $Q$, as well $\,\rg\in D(\beta_\Gamma)\,$ 
and $\,\xi_\Gamma\in\beta_\Gamma(\rg)\,$ almost
everywhere on $\Sigma$. In particular, thanks to {\bf (A5)}, the values of 
$\rho$ belong almost everywhere to the domain of definition of the original
function $g$. 
Hence, we may pass to the limit as $\,\eps_n\searrow0$
in the system \eqref{eps1}--\eqref{eps5} to infer that 
the quintuple $\,(\mu,\rho,\rg,\xi,\xi
_\Gamma)\,$ is a solution to the system \eqref{ss1}--\eqref{ss7}.

Finally, we have that $\pn\mu=0$ \gianni{a.e.}\ on $\Sigma$ and, in addition,
that $\,g(\rho),g'(\rho)\in C^0(\overline Q)$ and $\mu\in L^\infty(Q)$. Since
\beq\label{nuovo}
\pt\mu\,-\,(1+2g(\rho))^{-1}\Delta\mu\,=\,-(1+2g(\rho))^{-1}\,\mu\,g'(\rho)\,\pt\rho,
\eeq
we may view \eqref{ss1} as a linear uniformly parabolic equation for $\mu$ with continuous coefficients
and a \rhs~in $L^\infty(0,T;H)\cap L^2(0,T;L^6(\oma))$. As $\mu_0\in W$, it follows from optimal $L^p$-$L^q$-regularity
results (cf. \cite[Thm.~2.3]{Denk}) that $\,\mu\in W^{1,p}(0,T;H)\cap L^p(0,T;W)$, for all
$p\in [1,+\infty)$. We have thus shown the existence of a solution having the asserted regularity properties.

\vspace{2mm}\noindent
\pier{{}\underline{\sc Step 5:} \quad 
It remains to} prove the uniqueness result. For this purpose, we follow closely the proof of 
\gianni{\cite[Thm.~3.9]{CGPS5}} for the case that $\,\rho\,$ satisfies a zero Neumann condition, 
pointing out the differences originating from the dynamic boundary condition \eqref{ss5}. First, \pier{it is not difficult to check} that
Eq.~\eqref{ss1} can be rewritten as 
\beq\label{neu1}
\pt(\mu/\alpha(\rho))\,-\,\alpha(\rho)\,\Delta\mu\,=\,0 \quad\mbox{a.\,e. in }\,Q,
\eeq
where $\,\alpha:\overline{D(\beta)}\to (0,+\infty)\,$ is given by
\beq\label{defal}
\alpha(r)\,:=\,\bigl(1+2g(r) \bigr)^{-1/2}\quad\mbox{for \,$r\in \overline{D(\beta)}$}.
\eeq

Now suppose that two solutions $\,(\mu_i,\rho_i,\rho_{i_\Gamma},\xi_i,\xi_{i_\Gamma})$, $i=1,2$, to the
system \eqref{ss1}--\eqref{ss7} with the regularity \eqref{regmu}--\eqref{regxi} are given such that 
$\,\mu_i\ge 0\,$ almost everywhere in $Q$, for $i=1,2$. We put
$\,a_i:=\alpha(\rho_i)$, $z_i:=\mu_i/a_i$, for $i=1,2$, and set
\begin{align*}
&\mu:=\mu_1-\mu_2, \quad\rho:=\rho_1-\rho_2, \quad\rg:=\rho_{1_\Gamma}-\rho_{2_\Gamma},
\quad\xi:=\xi_1-\xi_2,\quad\xi_\Gamma:=\xi_{1_\Gamma}-\xi_{2_\Gamma},\\[1mm]
&a:=a_1-a_2,\quad z:=z_1-z_2,
\end{align*}
where we notice that $\,\rho_\Gamma=\rho_{|\Sigma}\,$ almost everywhere on $\Sigma$. Moreover, $z_i$
is bounded (since $\mu_i$ and $\rho_i$ are bounded), for $i=1,2$, and we have that 
\beq\label{L46}
|\nabla\rho_i|,|\nabla\mu_i|,|\nabla z_i|\in L^4(0,T;L^6(\oma)), \quad\mbox{for }\,i=1,2.
\eeq
In the following estimates, we denote by $\,C\,$ or $\,\widehat C\,$ positive constants that may
depend on the data of the system and on finite norms of the quantities $\mu_i,\rho_i,z_i$,
$i=1,2$. The constants $\,C\,$ may change their meaning within formulas and even within lines.  

At first, we write Eq.~\eqref{neu1} for $\mu_i,\rho_i$, $i=1,2$, and subtract the equations from each other,
which yields the identity $\,\,z_t-a_1\,\Delta(a_1z_1)+a_2\,\Delta(a_2z_2)=0$. \pier{Note that both $\,a_1z_1\,$ and $\,a_2z_2\,$ satisfy the Neumann homogeneous boundary condition \eqref{ss2}. Then, multiplication} by $\,z\,$
and integration over $Q_t$, where $t\in (0,T]$, leads to the equation
\beq\label{neu2}
\frac 12\xinto|z(t)|^2\dx\,+\,\txinto\bigl(\nabla(a_1z)\cdot\nabla(a_1z_1)-\nabla(a_2z)\cdot\nabla(a_2z_2)
\bigr)\dx\ds \gianni{{}=0}\,, 
\eeq 
whence, arguing exactly as between the formulas (6.8) and (6.10) in the proof of \gianni{\cite[Thm.~3.9]{CGPS5}},
we conclude the estimate
\begin{align}
\label{neu3}
&\frac 12\xinto|z(t)|^2\dx\,+\,\frac 12\txinto|\nabla(a_1z)|^2\dx\ds\,\le\,(1+\widehat C)\txinto|\nabla\rho|^2\dx\ds\non\\
&+\,C\txinto|\nabla\mu_2|^2\,|z|^2\dx\ds\,+\,C\txinto\left(|\nabla\rho_1|^2+|\nabla\rho_2|^2+|\nabla\mu_2|^2
+|\nabla z_2|^2\right)|\rho|^2\dx\ds\,.
\end{align}

Next, we write the equations \eqref{ss3} and \eqref{ss5} for the two solutions and subtract to obtain the identities
\begin{align}
\label{neu4}
&\rho_t-\Delta\rho+\xi=\pi(\rho_2)-\pi(\rho_1)+\mu\,g'(\rho_1)+\mu_2(g'(\rho_1)-g'(\rho_2))\quad\mbox{a.\,e. in }\,Q,
\\[1mm]
\label{neu5}
&\pn\rho+\pt\rg-\Delta_\Gamma\rg+\xi_\Gamma=\pig(\rho_{2_\Gamma})-\pig(\rho_{1_\Gamma})\quad\mbox{a.\,e. on }\,\Sigma.
\end{align}
Thus, multiplying \eqref{neu4} by $\rho$ and integrating over $Q_t$, where $t\in (0,T]$, we easily derive
from Young's inequality the estimate
\begin{align}
\label{neu6}
&\frac 12\left(\|\rho(t)\|_H^2+\|\rg(t)\|_{\Hg}^2\right)\,+\txinto|\nabla\rho|^2\dx\ds
\,+\tginto|\nabla_\Gamma\rg|^2\dg\ds\non\\
&+\txinto\xi\,\rho\dx\ds\,+\tginto\xi_\Gamma\,\rg\dg\ds\non\\
&\le\,C\txinto\left(|\mu|^2+|\rho|^2\right)\dx\ds\,+\,C\tginto|\rg|^2\dg\ds\,,
\end{align}
where, owing to the monotonicity of $\beta$ and $\beta_\Gamma$, the last two summands on the left-hand side
are nonnegative. At this point, we multiply the inequality \eqref{neu6} by $\,2+\widehat C\,$ and add the result to
\eqref{neu3}. Using the estimate
\beq\label{help!}
|\mu|=|a_1z_1-a_2z_2|\,\le\,|a||z_1|+|a_2||z|\,\le\,C\,(|\rho|+|z|),
\eeq
we then easily deduce that
\begin{align}
\label{neu7}
&\|z(t)\|_H^2+\|\rho(t)\|_H^2+\|\rg(t)\|^2_{\Hg}\non\\
&+\txinto\left(|\nabla(a_1z)|^2+|\nabla\rho|^2\right)\dx\ds
\,+\tginto|\nabla_\Gamma\rg|^2\dg\ds\non\\
&\le\, C\txinto\left(1+|\nabla\rho_1|+|\nabla\rho_2|+|\nabla\mu_2|+|\nabla z_2|\right)^2(|z|^2+|\rho|^2)\dx\ds \non\\
&\quad+\,C\tginto|\rg|^2\dg\ds\,.
\end{align}
Now, by virtue of \eqref{L46}, the function $\,k:=1+|\nabla\rho_1|+|\nabla\rho_2|+|\nabla\mu_2|
+|\nabla z_2|\,$ belongs to $L^4(0,T;L^6(\oma))$. A direct application of \gianni{\cite[Lem.~6.2]{CGPS5}} then
yields that the first summand on the right-hand side of \eqref{neu7} is bounded by an expression of the
form
\begin{align}
\label{neu8}
&\frac 12\txinto(|\nabla(a_1z)|^2+|\nabla\rho|^2)\dx\ds\non\\
&+\,C\tint\left(1+\|\nabla\rho_1(s)\|_6^4+\|k(s)\|_6^4\right)\left(\|z(s)\|_H^2+\|\rho(s)\|_H^2\right)\ds\,.
\end{align}
Combining this with \eqref{neu7}, we thus can infer from Gronwall's lemma that $z=\rho=0$ a.\,e. in $Q$ 
and $\rho_\Gamma=0$ a.\,e. on $\Sigma$. Hence, by \eqref{help!}, also $\mu=0$,
and, by comparison in \eqref{neu4}, $\xi=0$, a.\,e. in $Q$. Finally, by comparison in \eqref{neu5},
$\xi_\Gamma=0$ a.\,e. on $ \Sigma$. This concludes the proof of uniqueness. The assertion
is thus completely shown. 
 \qed

\vspace{5mm}\noindent
{\sc Remark 3.1:} \quad \,Since the expression on the right-hand side of \eqref{nuovo}
belongs to the space $L^2(0,T;L^6(\oma))$, it follows from the optimal $L^p$-$L^q$-regularity
result of \cite[Thm.~2.3]{Denk}) that also  
$\,\mu\in H^1(0,T;L^6(\oma))\cap L^2(0,T;W^{2,6}(\oma))$.

\section{Stability and further regularity}
\setcounter{equation}{0}

\noindent This section is devoted to the proofs of Theorem 2.3 and Theorem 2.4.

\vspace{2mm}\noindent
{\sc Proof of Theorem 2.3:}  \quad\,\,Assume that $\ug\in L^\infty(\Sigma)$, $\beta^\circ(\rho_0)\in L^\infty(\oma)$,
and $\beta^\circ(\rho_{0_\Gamma})\in L^\infty(\Gamma)$. We need to show that the unique solution $(\mu,\rho,
\rho_\Gamma,\xi,\xi_\Gamma)$
has the additional property that $\xi\in L^\infty(Q)$.
To this end, we recall that the approximating system \eqref{eps1}--\eqref{eps5}
introduced in Step~3 of the proof of Theorem 2.1 has a solution triple $\,(\mu^\eps,\rho^\eps,
\rho^\eps_\Gamma)$\, satisfying the estimate \eqref{K65}. Moreover,  
$\,\{\rho^\eps\}\,$ is bounded in $\,C^0(\overline{Q})\,$ uniformly with respect to $\eps>0$ , which entails, in particular, that
$\,\beta^\eps(\rho^\eps)\in C^0(\overline{Q})$, and thus $\,\beta^\eps(\rho^\eps_\Gamma)\in C^0(\overline{\Sigma})\,$
and $\,\beta_\Gamma^\eps(\rho^\eps_\Gamma) \in C^0(\overline{\Sigma})$. We therefore may
multiply \eqref{eps3} by $\,(\beta^\eps(\rho^\eps))^{2k-1}$, for every fixed $k\in\nz$. Integration over
$Q_t$, where $0<t\le T$, yields the identity  
\begin{align}\label{eq41}
&\xinto\Phi^\eps(x,t)\dx\,+\ginto\Phi_\Gamma^\eps(x,t)\dg\,+\,(2k-1)\txinto|\beta^\eps(\rho^\eps)|^{2k-2}
\,(\beta^\eps)'(\rho^\eps)|\nabla\rho^\eps|^2\dx\ds\non\\
&+(2k-1)\tginto|\beta^\eps(\rho^\eps_\Gamma)|^{2k-2}
\,(\beta^\eps)'(\rho^\eps_\Gamma)|\ng\rho^\eps_\Gamma|^2\dg\ds\,+\txinto|\beta^\eps(\rho^\eps)|^{2k}\dx\ds\non\\
&=\xinto\Phi^\eps(x,0)\dx\,+\ginto\Phi_\Gamma^\eps(x,0)\dg\,+\txinto(\beta^\eps(\rho^\eps))^{2k-1}\,
(\mu^\eps\,g'(\rho^\eps)-\pi(\rho^\eps))\dx\ds\non\\
&\quad +\tginto(\beta^\eps(\rho^\eps_\Gamma))^{2k-1}\,(-\beta_\Gamma^\eps(\rho^\eps_\Gamma)-\pig(\rho^\eps_\Gamma)
+\ug)\dg\ds,
\end{align}
with the functions
\beq
\gianni{\Phi^\eps}(x,t):=\int_0^{\rho^\eps(x,t)}(\beta^\eps(r))^{2k-1}\,{\rm d}r,\quad
\gianni{\Phi^\eps_\Gamma}(x,t):=\int_0^{\rho^\eps_\Gamma(x,t)}(\beta^\eps(r))^{2k-1}\,{\rm d}r.
\eeq
Since the functions $\,\beta^\eps\,$ and $\,\beta^\eps_\Gamma\,$ vanish at zero and are monotone, all of the
expressions on the \lhs~are nonnegative. We estimate the terms of the \rhs~individually, where we denote by
$\,C_i$, $i\in\nz$, positive constants that may depend on the data of the system but not on $\eps>0$. At first,
invoking \gianni{\bf(A1)}, \eqref{domeps}, and the assumption 
\pier{\eqref{hpxibdd}
that implies} $\beta^\circ(\rho_0)\in L^\infty(\oma)$, we obtain for the first summand on the \rhs~the estimate
\begin{align}\label{phi1}
&\xinto\!\int_0^{\rho_0(x)}(\beta^\eps(r))^{2k-1}\,{\rm d}r\dx
\,\le\,\xinto|\rho_0(x)|\,|\beta^\eps(\rho_0(x))|^{2k-1}\dx\non\\[1mm]
&\le\,\|\rho_0\|_\infty\xinto|\beta^\circ(\rho_0(x))
|^{2k-1}\dx
\,\le\, C_1\,\|\beta^\circ(\rho_0)\|_\infty^{2k-1}\,\le\,C_2\,C_3^{2k}\,.
\end{align}
By the same token, \pier{we have that}
\beq\label{phi2}
\int_\Gamma \Phi_\Gamma^\eps(x,0)\dx \,\le\, C_4\,C_5^{2k}\,.
\eeq
Next, recall that $\,\,\|\mu^\eps\,g'(\rho^\eps)-\pi(\rho^\eps)\|_{L^\infty(Q)}\,\le\,C_6$. Hence, using Young's
inequality \eqref{Young} with $\,p=\frac{2k}{2k-1}$, $q=2k$, and $\,\gamma=(\frac p2)^{1/p}$, we obtain for
the third integral on the \rhs, which we denote by $I_1$, the following estimate:
\begin{align}\label{phi3}
&I_1\,\le\,C_6\txinto|\beta^\eps(\rho^\eps)|^{2k-1}\dx\ds\,\le\,\frac 12\txinto|\beta^\eps(\rho^\eps)|^{2k}
\dx\ds\,+\,\frac 1{2k}\left(\frac {2k-1}{k}\right)^{2k-1}\!C_6^{2k}\non\\
&\le \,\frac 12\txinto|\beta^\eps(\rho^\eps)|^{2k}
\dx\ds\,+\,C_7\,C_8^{2k}\,.
\end{align} 
Finally, since $\ug\in L^\infty(\Sigma)$ by assumption, we obtain that 
$\,\,\|\ug-\pig(\rg)\|_{L^\infty(\Sigma)}\,\le\,C_9$. Moreover,
we have, using \eqref{domination2}, and invoking Young's inequality similarly as above,
\begin{align}
&-\,(\beta^\eps(\rho^\eps_\Gamma))^{2k-1}\,\beta_\Gamma^\eps(\rho^\eps_\Gamma)\,=\,-\,
|\beta^\eps(\rho^\eps_\Gamma)|^{2k-1}
\,|\beta^\eps_\Gamma(\rho^\eps_\Gamma)|\,\le\,-\frac 1\eta\,|\beta^\eps(\rho^\eps_\Gamma)|^{2k}\,+\,
\frac{C_\Gamma}\eta\,|\beta^\eps(\rho^\eps_\Gamma)|^{2k-1}\non\\
&\le\,-\,\frac 1{2\eta}\,|\beta^\eps(\rho^\eps_\Gamma)|^{2k}\,+\,C_{10}\,C_{11}^{2k}\,.
\end{align}
Hence, the last integral on the \rhs\ \gianni{of \eqref{eq41}}, which we denote by $I_2$, admits an estimate of the form
\beq\label{phi4}
I_2\,\le\,-\,\frac 1{2\eta}\tginto|\beta^\eps(\rho^\eps_\Gamma)|^{2k}\dg\ds\,+\,C_{12}\,C_{13}^{2k}\,.
\eeq
Combining the above estimates, we have thus shown that
\beq\label{phi5}
\|\beta^\eps(\rho^\eps)\|_{L^{2k}(Q)}^{2k}\,+\,\|\beta^\eps(\rho^\eps_\Gamma)\|_{L^{2k}(\Sigma)}^{2k}
\,\le\,C_{\gianni{14}}\,C_{\gianni{15}}^{2k} \quad\mbox{for all $k\in\nz$}.
\eeq
Taking the $(2k)$-th root on both sides of this inequality and then letting $k\to\infty$, we conclude
that
\beq\label{phi6}
\|\beta^\eps(\rho^\eps)\|_{L^\infty(Q)}\,+\,\|\beta^\eps(\rho^\eps_\Gamma)\|_{L^\infty(\Sigma)}\,\le\,\juerg{C_{16}}.
\eeq
But this means that the convergence result \eqref{coneps4} in Step~4 of the proof
of Theorem 2.1 can be replaced by the stronger statement
\begin{align}
\beta^{\eps_n}(\rho^{\eps_n})\to \xi \quad\mbox{weakly-star in}\,L^\infty(Q),
\end{align}
so that the solution constructed there satisfies $\,\xi\in L^\infty(Q)$. Since the solution is
unique, the assertion is proved.\qed

\vspace{5mm}
Before proving Theorem 2.4, we cite a known auxiliary result (cf. \cite[Thm.~2.2]{CS}) that will be used
during the course of the proof.

\vspace{4mm}\noindent
{\sc Lemma 4.1:} \quad\,\,{\em Suppose that functions $a\in L^\infty(Q)$,
$a_\Gamma\in L^\infty(\Sigma)$, $\sigma\in L^2(Q)$, and $\sigma_\Gamma\in L^2(\Sigma)$ are 
given. Then the linear initial-boundary value problem}
\begin{align}
\label{CoSp1}
&\pt y-\Delta y\,+\,a\,y\,=\,\sigma \quad\mbox{{\em a.\,e. in }}\,Q,\\
\label{CoSp2}
&\pn y+\pt y_\Gamma-\Delta_\Gamma y_\Gamma\,+\,a_\Gamma\,y_\Gamma\,=\,\sigma_\Gamma,\quad y_\Gamma
=y_{\pier{|\Sigma}}, \quad\mbox{{\em a.\,e. on }}\,\Sigma,\\
\label{CoSp3}
&y(0)=0 \quad\mbox{{\em a.\,e. in }}\,\oma,\quad y_\Gamma(0)=0
\quad\mbox{{\em a.\,e. on }}\,\Gamma, 
\end{align}
{\em has a unique solution pair satisfying $\,y\in H^1(0,T;H)\cap C^0([0,T];V)\cap L^2(0,T;H^2(\oma))\,$ and
$y_\Gamma\in H^1(0,T;\Hg)\cap C^0([0,T];\Vg)\cap L^2(0,T;H^2(\Gamma))$. Moreover, there is a constant $C_L>0$,
which depends only on $\Omega$, $T$, $\|a\|_{L^\infty(Q)}$, and $\|a_\Gamma\|_{L^\infty(\Sigma)}$,  such that, for every $t\in (0,T]$,}
\begin{align}
\label{CoSp4}
&\|y\|_{H^1(0,t;H)\cap C^0([0,t];V)\cap L^2(0,t;H^2(\oma))}
\,+\,\|y_\Gamma\|_{H^1(0,t;\Hg)\cap C^0([0,t];\Vg)\cap L^2(0,t;H^2(\Gamma))}\non\\
&\le\,C_L\left(\|\sigma\|_{L^2(Q_t)}\,+\,\|\sigma_\Gamma\|_{L^2(\Sigma_t)}\right)\,.
\end{align}

\vspace{3mm}
We are now prepared to prove Theorem 2.4.

\vspace{3mm}\noindent
{\sc Proof of Theorem 2.4:} \quad
\gianni{We first prove~\eqref{separation}
by recalling that both $\rho$ and $\xi$ are bounded
(see Rem.~2.2 and Thm.~2.3).
We construct~$r^*$.
If $r_+=+\infty$ we can take $r^*:=\sup\rho<r_+$.
If instead $r_+<+\infty$, we choose $s^*\geq0$ such that $s^*\geq\xi$ a.e.\ in~$Q$.
By observing that $\lim_{r\nearrow r_+}\beta^\circ(r)=+\infty$ by maximal monotonicity,
we find $r^*\in(0,r_+)$ such that $\beta^\circ(r)>s^*$ for every $r\in(r^*,r_+)$.
Then, it holds that $\rho\leq r^*$ on~$\overline Q$.
Indeed, if the continuous function $\rho$ satisfies $\rho(x,t)=r^*+2\delta$ for some $\delta>0$ at some point~$(x,t)$, 
we should have $\rho\geq r^*+\delta$ in some subset $Q'\subset Q$ with positive measure.
This implies that $\xi\geq\beta^\circ(r^*+\delta)>s^*$ a.e.\ in~$Q'$,
and this contradicts the definition of~$s^*$.
As $r_*$ can be obtained similarly, the separation property \eqref{separation} is established.
From this and $D(\beta_\Gamma)=(r_-,r_+)$ it immediately follows that even $\xi_\Gamma$ is bounded.}

We \gianni{assume that the potentials are smooth and} show the validity of \eqref{stabu}. To this end,
suppose that \gianni{{$\juerg{u_{i_\Gamma}}\in H^1(0,T;\Hg)\cap L^\infty(\Sigma)$}, $i=1,2$, are given, and
let $(\mu_i,\rho_i,\rho_{i_\Gamma}, \xi_i,\xi_{i_\Gamma})$, where $\xi_i=\beta(\rho_i)$ and
$\xi_{i_\Gamma}=\beta_\Gamma(\rho_{i_\Gamma})$, be corresponding solutions to \eqref{ss1}--\eqref{ss7},
for $i=1,2$, that enjoy the regularity properties \eqref{regmu}--\eqref{regxi}. 
In the following, we denote by $C>0$ constants that depend only on
the data of the system, on $\juerg{\|u_{i_\Gamma}\|}_{H^1(0,T;\Hg)}$, $i=1,2$, and on the 
norms of $\,(\mu_i,\rho_i,\rho_{i_\Gamma})$, $i=1,2$, in the spaces indicated in \eqref{regmu}--\eqref{regrg}. Now, we put
\begin{align*}
\mu:=\mu_1-\mu_2,\quad\rho:=\rho_1-\rho_2, \quad\rg:=\rho_{1_\Gamma}-\rho_{2_\Gamma},
\quad\ug:=\juerg{u_{1_\Gamma}-u_{2_\Gamma}}.   
\end{align*}
Obviously, $\rg=\rho_{|\Sigma}$ on $\Sigma$. \gianni{Moreover, both $\rho_i$ and $\rho_{i_\Gamma}$
satisfy \eqref{separation} for $i=1,2$.
As all the functions $g$, $g'$, $\pi$, $\beta$, $\pi_\Gamma$ and $\beta_\Gamma$ 
are Lipschitz continuous on~$[r_*,r^*]$, 
we have that}
\begin{align}
\label{b1}
&\max_{0\le i\le 1}\left|g^{(i)}(\rho_1)-g^{(i)}(\rho_2)\right|+
\left|\pi(\rho_1)-\pi(\rho_2)\right|
+\left|\beta(\rho_1)-\beta(\rho_2)\right|\le\,C\,|\rho| \quad\mbox{in }\,\overline Q,\\[1mm]
\label{b2}
&\left|\beta_\Gamma(\rho_{1_\Gamma})-\beta_\Gamma(\rho_{2_\Gamma})\right|+\left
|\pig(\rho_{1_\Gamma})-\pig(\rho_{2_\Gamma})\right|\,\le\,C\,|\rg|\quad\mbox{on \,}
\overline\Sigma.
\end{align}
\gianni{Furthermore}, we easily verify that 
\begin{align}
\label{diff1}
&(1+2g(\rho_1))\,\pt\mu\,+\,g'(\rho_1)\,\pt\rho_1\,\mu\,-\,\Delta\,=\,-\pier{{}2{}}\,\pt\mu_2\,
(g(\rho_1)-g(\rho_2))\mu\non\\
&\hspace*{30mm}-\,\mu_2\,(g'(\rho_1)-g'(\rho_2))\,\pt\rho_1\,-\,\mu_2\,g'(\rho_2)\,\pt\rho
\quad\mbox{a.\,e. in $\,Q$},\\[1mm]
\label{diff2}
&\pn\mu\,=\,0\quad\mbox{a.\,e. on }\,\Sigma,\quad \mu(0)=0\quad\mbox{a.\,e. in }\,\oma,\\[1mm]
\label{diff3}
&\pt\rho-\Delta\rho\,=\,\mu_1\,(g'(\rho_1)-g'(\rho_2))+g'(\rho_2)\,\mu + \pi(\rho_2)-\pi(\rho_1)
\non\\
&\hspace*{30mm}+\,\beta(\rho_2)-\beta(\rho_1)\quad
\mbox{a.\,e. in $\,Q$},\\[1mm]
\label{diff4}
&\pn\rho+\pt\rg-\Dg\rg\,=\,\beta_\Gamma(\rho_{2_\Gamma})-\beta_\Gamma(\rho_{1_\Gamma})+
\pig(\rho_{2_\Gamma})-\pig(\rho_{1_\Gamma})+\ug\quad\mbox{a.\,e. on $\,\Sigma$},\\[1mm]
\label{diff5}
&\rho(0)=0\quad\mbox{a.\,e. in $\,\oma$},\quad\rho_\Gamma(0)=0\quad\mbox{a.\,e. on }\,\Gamma.
\end{align}

\vspace{2mm}\noindent
First, we observe that $\,\,\pt((\frac 12+g(\rho_1))\,\mu^2)=(1+2g(\rho_1))\,\mu\,\pt\mu+
g'(\rho_1)\,\pt\rho_1\,\mu^2$. Hence, if we multiply \eqref{diff1} by $\,\mu\,$ and
integrate over $Q_t$, where $t\in (0,T]$, then we obtain the estimate
\begin{align}
\label{est41}
\xinto\left(\mbox{$\frac 12$}+g(\rho_1(t))\right)|\mu(t)|^2\dx
\,+\txinto|\nabla\mu|^2\dx\ds\,\le\,I_1+I_2+I_3,
\end{align}
where the expressions $I_j$, $1\le j\le 3$, are defined and estimated as follows: at first,
we use H\"older's and Young's inequalities, the continuity of the embedding
$V\subset L^4(\oma)$, and \eqref{b1}, to obtain that
\begin{align}
\label{est42}
I_1&=\,-\txinto \pier{{}2{}}\,\pt\mu_2\,(g(\rho_1)-g(\rho_2))\,\mu\dx\ds
\,\le\,C\tint\|\pt\mu_2(s)\|_2\,\|\rho(s)\|_4\,\|\mu(s)\|_4\ds\non\\
&\le\,\frac 12\tint\|\mu(s)\|_V^2\ds\,+\,C\tint\|\pt\mu_2(s)\|_H^2\,\|\rho(s)\|_V^2\ds\,.
\end{align}
By the same token, \pier{we infer that}
\begin{align}
\label{est43}
I_2&=\,-\txinto \mu_2\,(g'(\rho_1)-g'(\rho_2))\,\pt\rho_1\,\mu\dx\ds
\,\le\,C\tint\|\pt\rho_1(s)\|_4\,\|\rho(s)\|_4\,\|\mu(s)\|_2\ds\non\\
&\le\,C\txinto|\mu|^2\dx\ds\,+\,C\tint\|\pt\rho_1(s)\|_V^2\,\|\rho(s)\|_V^2\ds\,.
\end{align}
Finally, \pier{we have that}
\beq\label{est44}
I_3=\,-\txinto\mu_2\,g'(\rho_2)\,\pt\rho\,\mu\dx\ds\,\le\,\frac 14\txinto|\pt\rho|^2\dx\ds\,
+\,C\txinto|\mu|^2\dx\ds\,.
\eeq
Combining \eqref{est41}--\eqref{est44}, and recalling that $g$ is nonnegative, we have thus
shown an estimate of the form
\begin{align}
\label{est45}
&\frac 12\xinto|\mu(t)|^2\dx\,+\,\frac 12\tint\|\mu(s)\|_V^2\ds 
\,\le\,\frac 14\txinto|\pt\rho|^2\dx\ds\non\\
&\quad +\,C\tint\left(1+\|\pt\mu_2(s)\|_H^2\,+\,\|\pt\rho_1(s)\|_V^2\right)
\left(\|\rho(s)\|_H^2\,+\,\|\mu(s)\|_H^2\right)\ds\,.
\end{align}

Now, we add $\,\rho\,$ to both sides of \eqref{diff3} and $\,\rg\,$ to
both sides of \eqref{diff4}. Then we multiply \eqref{diff1} by $\pt\rho$ and 
integrate over $Q_t$, where $0<t\le T$. Using \eqref{b1} and \eqref{b2}, we deduce
that
\begin{align}\label{est46}
&\txinto|\pt\rho|^2\dx\ds\,+\,\tginto|\pt\rg|^2\dg\ds\,+\,\frac 12\,
\bigl(\|\rho(t)\|_V^2 \,+\,\|\rg(s)\|_{\Vg}^2\bigr)\non\\
&\le\,C\txinto|\pt\rho|\bigl(|\mu|\,+\,|\rho|\bigr)\dx\ds\,+\,C\tginto
|\pt\rg|\,\bigl(|\rg|\,+\,|\ug|\bigr)\dg\ds\non\\
&\le\,\frac 14\txinto|\pt\rho|^2\dx\ds\,+\,\frac 12\tginto|\rg|^2\dx\ds
\,+\,C\tginto|\ug|^2\dg\ds\non\\
&\quad\, +\,C\txinto(|\rho|^2\,+|\mu|^2)\dx\ds\,+\,C\tginto|\rg|^2\dg\ds\,.
\end{align}
Now\pier{, we add \eqref{est45} to \juerg{the} final result of \eqref{est46}. O}bserve that the mapping \,\,$s \mapsto \|\pt\mu_2(s)\|_H^2\,+\,\|\pt\rho_1(s)\|_V^2$
is known to belong to $L^1(0,T)$. Hence, invoking Gronwall's lemma, we can infer
from \eqref{est45} and \eqref{est46} the estimate
\begin{align}\label{est47}
&\|\mu\|_{L^\infty(0,t;H)\cap L^2(0,t;V)}+\|\rho\|_{H^1(0,t;H)\cap L^\infty(0,t;V)} \non \\
&\quad +\|\rg\|_{H^1(0,t;\Hg)\cap L^\infty(0,t;\Vg)}\,\le\,C\,\|\ug\|_{L^2(0,t;\Hg)}\,.
\end{align}

At this point, we observe that the system \eqref{diff3}--\eqref{diff5} is of the form
\eqref{CoSp1}--\eqref{CoSp3}, \pier{with} $y:=\rho$, $y_\Gamma:=\rg$, $a:=1$, $a_\Gamma:=1$,
and where $\sigma$ and $\sigma_\Gamma$ are given by the right-hand sides of
\eqref{diff3} and \eqref{diff4}, respectively. Therefore, invoking \eqref{b1},
\eqref{b2}, and \eqref{est47}, we conclude from Lemma 4.1 that
\begin{align}
\label{est48}
&\|\rho\|_{H^1(0,t;H)\cap C^0([0,t];V)\cap L^2(0,t;H^2(\oma))}
\,+\,\|\rg\|_{H^1(0,t;\Hg)\cap C^0([0,t];\Vg)\cap L^2(0,t;H^2(\Gamma))}\non\\
&\le\,C\,\|\ug\|_{L^2(0,t;\Hg)}\,.
\end{align} 
With this, the stability estimate \eqref{stabu} is shown.\qed

\pier{\small%
 
}%
\end{document}